\documentclass{amsart}

\newtheorem{theorem}{Theorem}[section]
\newtheorem*{teorema}{Theorem}
\newtheorem{corollary}[theorem]{Corollary}
\newtheorem{lemma}[theorem]{Lemma}
\newtheorem{proposition}[theorem]{Proposition}

\newtheorem*{claim}{Claim}
\numberwithin{equation}{section}
\newcommand{\N}{\mathbb N}

\newcommand{\Q}{\mathbb Q}
\newcommand{\C}{\mathbb C}
\newcommand{\F}{\mathbb F}
\newcommand{\A}{\mathbb A}
\newcommand{\K}{{\sf K}}

\newcommand{\fq}{\F_{\hskip-0.7mm q}}
\newcommand{\cfq}{\overline{\F}_{\hskip-0.7mm q}}

\newcommand{\slp} {straight--line program}
\newcommand{\slps} {straight--line programs}

\newcommand{\lfiberr} {\ifm {V_{P^{(r)}}}{$V_{P^{(r)}}$}}
\newcommand{\lfibers}  {\ifm {V_{P^{(s)}}}{$V_{P^{(s)}}$}}
\newcommand{\lfibersstar}  {\ifm {V_{P^{(s)}}^*}{$V_{P^{(s)}}^*$}}
\newcommand{\lfiberspone} {\ifm {V_{P^{(s+1)}}}{$V_{P^{(s+1)}}$}}
\newcommand{\lfibersponestar} {\ifm {V_{P^{(s+1)}}^*}{$V_{P^{(s+1)}}^*$}}
\newcommand{\lcurvespone} {\ifm {W_{P^{(s+1)}}}{$W_{P^{(s+1)}}$}}
\newcommand{\lcurvesponestar} {\ifm {W_{P^{(s+1)}}^*}{$W_{P^{(s+1)}}^*$}}

\newcommand{\LLambda}{\ifm {\mathcal{L}_\Lambda}{$\mathcal{L}_\Lambda$}}
\newcommand{\Llambda}{\ifm {\mathcal{L}_\lambda}{$\mathcal{L}_\lambda$}}
\newcommand{\llambda}{\ifm {\ell_\lambda}{$\ell_\lambda$}}
\newcommand{\Llambdaone}{\ifm {\mathcal{L}_{\lambda_1}}{$\mathcal{L}_{\lambda_1}$}}
\newcommand{\Llambdatwo}{\ifm {\mathcal{L}_{\lambda_2}}{$\mathcal{L}_{\lambda_2}$}}

\def\ifm#1#2{\relax \ifmmode#1\else#2\fi}

\newcommand{\klk}    {\ifm {,\ldots,} {$,\ldots,$}}

\newcommand{\plp}    {\ifm {+\cdots+} {$+\ldots+$}}

\newcommand{\om}[2]   {{#1}_1 \klk {#1}_{#2}}

\newcommand{\sspar} {\vskip 0.125cm}
\newcommand{\spar} {\vskip 0.25cm}

\newcommand{\xo}[1]  {\ifm {\om X {#1}} {$\om X {#1}$}}
\newcommand{\xon}    {\ifm {\om X n} {$\om X n$}}
\newcommand{\yo}[1]  {\ifm {\om Y {#1}} {$\om Y {#1}$}}
\newcommand{\yon}    {\ifm {\om Y n} {$\om Y n$}}
\newcommand{\fo}[1]  {\ifm {\om F {#1}} {$\om F {#1}$}}
\begin{document}

\title[Computation of a rational point]{Fast computation of a
rational point of a variety over a finite field}%
\author{A. Cafure${}^{1}$}
\address{${}^{1}$Departamento de Matem\'atica, Facultad de
Ciencias Exactas y Naturales, Universidad de Buenos Aires, Ciudad
Universitaria, Pabell\'on I (1428) Buenos Aires, Argentina.}
\email{acafure@dm.uba.ar}

\author{G. Matera${}^{2,3}$}
\address{${}^{2}$Instituto de Desarrollo Humano,
Universidad Nacional de General Sarmiento, J.M. Guti\'errez 1150
(1613) Los Polvorines, Buenos Aires, Argentina.}
\email{gmatera@ungs.edu.ar}
\address{${}^{3}$Member of 
the CONICET, Argentina.}

\thanks{Research was partially supported by the
following grants: UBACyT X198, PIP CONICET 2461 and UNGS 30/3005}

\date{\today}%
\begin{abstract}
We exhibit a probabilistic algorithm which computes a rational
point of an absolutely irreducible variety over a finite field
defined by a reduced regular sequence. Its time--space complexity
is roughly quadratic in the logarithm of the cardinality of the
field and a geometric invariant of the input system (called its
degree), which is always bounded by the B\'ezout number of the
system. Our algorithm works for fields of any characteristic, but
requires the cardinality of the field to be greater than a
quantity which is roughly the fourth power of the degree of the
input variety.
\end{abstract}
\maketitle
%
%
\section{Introduction}
Let $p$ be a prime number, let $q:=p^k$, let $\fq$ be the finite
field of $q$ elements and let $\cfq$ denote its algebraic closure.
For a given $n\in\N$, let $\A^n$ denote the $n$--dimensional
affine space $\cfq^n$ endowed with its Zariski topology. Let be
given a finite set of polynomials
$F_1,\dots,F_m\in\fq[X_1,\dots,X_n]$ and let $V$ denote the affine
subvariety of $\A^n$ defined by $F_1,\dots,F_m$. In this paper we
consider the problem of computing a $q$--rational point of the
variety $V$, i.e. a point $x\in\fq^n$ such that $F_i(x)=0$ holds
for $i=1\klk m$.

This is an important problem of mathematics and computer science,
with many applications. It is NP--complete, even if the equations
are quadratic and the field considered is $\F_2$. Furthermore,
\cite{GaKaSh97} shows that determining the number of rational
points of a sparse plane curve over a finite field is
$\#$P--complete. In fact, several multivariate cryptographic
schemes based on the hardness of solving polynomial equations over
a finite field have been proposed and cryptoanalyzed (see e.g.
\cite{CoKlPaSh00}). The problem is also a critical point in areas
such as coding theory (see e.g. \cite{BoPe99}, \cite{LiNi83}),
combinatorics \cite{LiPi84}, etc.

There is not much literature on the subject. In \cite{GaShSi03},
an algorithm computing the set of $q$--rational points of a plane
curve over a finite field has been proposed. On the other hand,
\cite{KiSh99} and \cite{CoKlPaSh00} exhibit algorithms which
solve an overdetermined system of quadratic equations over a
finite field, based on a technique of linearization.

Algorithms finding rational points on a general variety over a
finite field are usually based on rewriting techniques (see e.g.
\cite{CoLiOS92}, \cite{CoLiOS98}). Unfortunately, such algorithms
have superexponential complexity, which makes them infeasible for
realistically sized problems. Indeed, their most efficient
variants (see e.g. \cite{Faugere02}) have worst--case complexity
higher than exhaustive search in polynomial equation systems over
$\F_2$ \cite{CoKlPaSh00}.

A different approach is taken in \cite{HuWo99}, which exhibits an
algorithm solving polynomial systems over a finite field by means
of deformations, based on a perturbation of the original system
and a subsequent path--following method. Nevertheless, the
perturbation typically introduces spurious solutions which may be
computationally expensive to identify and eliminate in order to
obtain the actual solutions. Furthermore, the algorithm is
algebraically robust or universal in the sense of
\cite{HeMaPaWa98} and \cite{CaGiHeMaPa03}, which implies
exponential lower bounds on its time complexity.

In the case of polynomial equation solving over the complex or
real numbers, the series of papers \cite{GiHeMoPa95},
\cite{Pardo95}, \cite{GiHeMoMoPa98}, \cite{GiHaHeMoMoPa97},
\cite{GiHeMoPa97}, \cite{BaGiHeMb97}, \cite{BaGiHeMb01} (see also
\cite{HeMaWa01}, \cite{GiLeSa01}, \cite{Lecerf03}) introduces a
new symbolic elimination algorithm. Its complexity is roughly the
product of the complexity of the input polynomials and a {\em
polynomial} function of a certain geometric invariant of the input
system, called its {\em degree}. Let us observe that the degree is
always bounded by the B{\'e}zout--number of the input system and
happens often to be considerably smaller.

Our intention is to develop a new family of elimination algorithms
for polynomial equation solving over finite fields of the type
above. For this purpose, and as a first step towards this aim, in
this article we exhibit a new probabilistic algorithm which solves
a critical problem of effective elimination over finite fields:
the computation of a rational point of an $\fq$--definable
absolutely irreducible variety. Our main result is summarized in
the following theorem (see Corollary \ref{coro:final} for a
precise complexity statement):
\begin{teorema}
Let $n\ge 3$ and $d\ge 2$. Let $\fo r\in\fq[\xon]$ be polynomials
of maximum degree $d$ which form a regular sequence of
$\fq[\xon]$. Suppose that $\fo s$ generate a radical ideal of
$\fq[\xon]$ for $1\le s\le r$ and let $V_s:=V(\fo s)\subset\A^n$.
Let $\delta:=\max_{1\le s\le r}\deg V_s$. Let be given a \slp\ in
$\fq[\xon]$ using space $\mathcal{S}$ and time $\mathcal{T}$ which
represents $\fo r$. Suppose further that $V:=V_r$ is absolutely
irreducible and $q>8n^2d\delta_r^4$ holds. Then there exists a
probabilistic algorithm which computes a $q$--rational point of
$V$ with space $O{\,\widetilde{\,}}\big(\mathcal{S}\delta^2\log^2
q\big)$ and time
$O{\,\widetilde{\,}}(\mathcal{T}\delta^2\log^2q)$. \spar

\rm (Here $O{\,\widetilde{\,}}\,$ refers to the standard Soft--Oh
notation which does not take into account logarithmic terms.
Further, we have ignored terms depending on $n$ and $d$.)
\end{teorema}
The complexity estimate of our algorithm is polynomial in the
degree of the system $\delta$ mentioned above and the logarithm of
 $q$.
Therefore, taking into account the {\em worst--case} estimate
$\delta\le D:=\prod_{i=1}^r\deg(F_i)$, we conclude that our
algorithm achieves, for the first time in polynomial equation
solving over finite fields, a complexity {\em polynomial} in the
B\'ezout number $D$ and $\log q$. In particular, our complexity
result exponentially improves the $d^{O(n^2)}\log^{O(1)}\!q$
worst--case estimates of \cite{HuWo99} and Gr\"obner solving
algorithms.

In the above statement we assume that the input polynomials $\fo
r$ form a {\em reduced} regular sequence, i.e. $\fo s$ generate a
radical ideal for $1\le s\le r$. Let us remark that this
hypothesis can be easily recovered from a regular sequence
generating a radical ideal by a generic linear combination of the
input polynomials (see e.g. \cite[Proposition 37]{KrPa96}).
Furthermore, using techniques inspired in \cite{Lecerf02},
\cite{Lecerf03} it is possible to extend our algorithm to {\em
arbitrary} polynomial equation systems over $\fq$ defining an
absolutely irreducible variety (this extension shall be considered
in a forthcoming work). \spar

Our algorithm may be divided into three main parts. The first part
is a procedure which has as input a reduced regular sequence $\fo{r}
\in\fq[\xon]$ and outputs a complete description of the generic
point of the input variety $V:=V(\fo{r})$. Such a description is
provided by a $\K$--definable generic linear projection
$\pi_r:V\to\A^{n-r}$ and a parametrization of an
unramified generic fiber $\pi_r^{-1}(P^{(r)})$, where $\K$ is a
suitable finite field extension of $\fq$ (cf. Sections
\ref{subsection:geo_sol}, \ref{subsection:lifting_points}).

In Section \ref{section:computation_geo_sol} we describe
this (recursive) procedure, which proceeds in $r-1$ steps.
Its $s$--th step computes a complete description of the
generic point of $V_{s+1}:=V(\fo{s+1})$, which is represented
by an unramified fiber $\pi_{s+1}^{-1}(P^{(s+1)})$ of a finite
$\K$--definable linear projection $\pi_{s+1}:V_{s+1}\to\A^{n-s-1}$.
For this purpose, in
Section \ref{subsec:lifting_to_curve} the unramified fiber
$\pi_{s}^{-1}(P^{(s)})$ of the previous step is ``lifted" to a
suitable curve $W_{P^{(s+1)}}$, contained in $V_s:=V(\fo{s})$, whose
intersection with the hypersurface defined by $F_{s+1}$ yields a
complete description of the fiber
$\pi_{s+1}^{-1}(P^{(s+1)})$. This intersection is computed in Sections
\ref{subsec:intersection_step} and \ref{subsec:shape_lemma}.

In the second part of our algorithm (Section
\ref{subsec:fq_def_geo_sol}), we obtain an $\fq$--definable
description of the generic point of $V$. For this purpose, we
develop a symbolic homotopy algorithm, based on a global
Newton--Hensel lifting, which ``moves" the $\K$--definable finite
morphism $\pi_{r}:V_{r}\to\A^{n-r}$ and the $\K$--definable
generic unramified fiber  $\pi_{r}^{-1}(P^{(r)})$ previously
obtained, into an $\fq$--definable finite morphism
$\pi:V\to\A^{n-r}$ and an  $\fq$--definable generic unramified
fiber $\pi^{-1}(P)$. Combining this procedure with an effective
version of the first Bertini theorem, in the third part of our
algorithm we obtain an absolutely irreducible plane $\fq$--curve
$\mathcal{C}$ with the property that any $q$--rational smooth
point of $\mathcal{C}$ immediately yields a $q$--rational point of
the input variety $V$ (see Section
\ref{section:comput_rat_point}). Then, in Section
\ref{subsec:comput_point_curve} we compute a $q$--rational point
of the curve $\mathcal{C}$ with a probabilistic algorithm which
combines the classical Weil's estimate and a procedure based on
factorization and gcd computations.\spar

A critical point of our algorithm is the {\em a priori}
determination of the linear projections $\pi_s$ and the points
$P^{(s)}$ for $1\le s\le r$. In Section \ref{section:preparation}
we show that this data can be generically chosen and obtain
explicit estimates on the degrees of the polynomials underlying
this genericity condition (improving significantly previous
estimates). Therefore, using the Zippel--Schwartz test
(see \cite{Zippel79}, \cite{Schwartz80} and Section
\ref{subsection:slps}) we may randomly find such
linear projections and points with high probability of success.
\spar

Let us remark that our algorithm does not impose any restriction
on the characteristic $p>0$, but requires the cardinality $q$ of
the field $\fq$  to satisfy the condition $q>8n^2d\delta_r^4$,
where $\delta_r$ is the degree of the variety $V$. Nevertheless,
it is clear that our algorithm cannot work unless there exists at
least one rational point of the input variety $V$. Since the
existence of a rational point of an absolutely irreducible variety
over $\fq$ of degree $\delta_r$ cannot be asserted (up to the
authors's knowledge) unless $q>2\delta_r^4$ holds, we see that our
condition on $q$ comes quite close to this ``minimal" requirement.

Finally, we observe that our algorithm can be efficiently extended
to the case of an $\fq$--definable variety $V$ with an absolutely
irreducible $\fq$--definable component of dimension equal to $\dim
V$. On the other hand, extensions to the general case of an
arbitrary variety over $\fq$ are likely to produce a significant
increase of the time--space complexity of our algorithm (see
\cite{HuWo99}).
%
%
\section{Notions and notations.}
\label{section:notions} We use standard notions and notations of
commutative algebra and algebraic geometry as can be found in e.g.
\cite{Kunz85}, \cite{Shafarevich84}, \cite{Matsumura80}.

Let $\fq$ and $\cfq$ denote the finite field of $q$ elements and its
algebraic closure respectively, and let $\K$ be a subfield
of $\cfq$ containing $\fq$. Let $\xon$ be indeterminates over $\K$
and let ${\K}[\xon]$ denote the ring of $n$--variate polynomials in
the indeterminates $\xon$ and coefficients in $\K$. Let $V$ be a
$\K$--definable affine subvariety of $\A^n$ (a $\K$--variety for
short). We shall denote by $I(V)\subset {\K}[\xon]$ its defining
ideal and by ${\K}[V]$ its coordinate ring, namely, the quotient
ring ${\K}[V]:={\K}[\xon]/I(V)$.

If $V$ is irreducible as a $\K$--variety ($\K$--irreducible for
short), we define its {\em degree} as the maximum number of points
lying in the intersection of $V$ with an affine linear subspace
$L$ of $\A^n$ of codimension $\dim(V)$ for which $\#(V\cap
L)<\infty$ holds. More generally, if $V=C_1\cup\cdots\cup C_h$ is
the decomposition of $V$ into irreducible $\K$--components, we
define the degree of $V$ as $\deg(V):=\sum_{i=1}^h\deg(C_i)$ (cf.
\cite{Heintz83}). In the sequel we shall make use of the following
{\em B\'ezout inequality} (\cite{Heintz83}; see also
\cite{Fulton84}): if $V$ and $W$ are $\K$--subvarieties of $\A^n$,
then
\begin{equation}\label{equation:Bezout}\deg (V\cap W)\le \deg V
\deg W.\end{equation}

A $\K$--variety $V \subset\A^n$ is {\em absolutely irreducible} if
it is irreducible as $\cfq$--variety.
%
%
\subsection{Geometric solutions}\label{subsection:geo_sol}

In order to describe the geometric aspect of our procedure we need
some more terminology, essentially borrowed from
\cite{GiHaHeMoMoPa97}. Let  us consider an equidimensional
$\K$--variety $W \subset \A^n$ of dimension $m\geq 0$ and degree
$\deg W$, defined by polynomials $F_1\klk F_{n-m} \in {\K}[\xon]$
which form a regular sequence. A {\em geometric solution} of $W$
consists of the following items:
\begin{itemize}
\item
a linear change of variables, transforming the variables
$X_1,\ldots,X_n$ into new ones, say $Y_1,\ldots,Y_n$, with the
following properties:
\begin{itemize}
\item the linear map $\pi:W\to\A^{m}$ defined by
$Y_1,\ldots,Y_{m}$ is a finite surjective morphism. In this case,
the change of variables is called a {\em Noether normalization} of
$W$ and we say that the variables $Y_1,\ldots,Y_n$ are in {\em
Noether position} with respect to $W$, the variables
$Y_1,\ldots,Y_{m}$ being {\em free}. The given Noether
normalization induces an integral ring extension
$R_m:=\cfq[Y_1,\ldots,Y_{m}]\hookrightarrow \cfq[W]$. Observe that
$\cfq[W]$ is a free $R_m$--module whose rank we denote by
$\mathrm{rank}_{R_m}\cfq[W]$. Notice that $\mathrm{rank}_{R_m}
\cfq[W]\le \deg W$ (see e.g. \cite{GiHeSa93}) and
$\cfq[W]\cong\cfq[X_1,\ldots,X_n]/(F_1,\ldots,F_{m-n})$ hold.
\item
the linear form $Y_{m+1}$ induces a primitive element of the ring
extension $R_m\hookrightarrow\cfq[W]$, i.e. an element
$y_{m+1}\in\cfq[W]$ whose (monic) minimal polynomial $q^{(m)}\in
R_m[T]$ over $R_m$ satisfies the condition
$\deg_Tq^{(m)}=\mathrm{rank}_{R_m}\cfq[W]$. Observe that we always
have $\deg q^{(m)}=\deg_Tq^{(m)}\le\deg W$.
\end{itemize}
\item the minimal polynomial $q^{(m)}$ of $y_{m+1}$ over $R_m$.
\item a generic ``{\em parametrization}" of the variety $W$ by the
zeroes of $q^{(m)}$, of the form $(\partial q^{(m)}/\partial T)(T)
Y_{m+2}-v_{m+2}^{(m)}(T) \klk (\partial q^{(m)}/\partial T)(T)
Y_n-v_n^{(m)}(T) $ with $v_{m+2}^{(m)}\klk v_n^{(m)}\in R_m[T]$.
We require that $ \deg_T\ v_{m+j}^{(m)}< \deg_T(q^{(m)})$ and $
(\partial q^{(m)}/\partial
T)(Y_{m+1})Y_{m+j}-v_{m+j}^{(m)}(Y_{m+1})\in (F_1\klk F_{n-m})$
hold for $2\le j\le n-m$. Observe that this parametrization is
unique up to scaling by nonzero elements of $\cfq$.
\end{itemize}

Observe that in the case that $W$ is a zero--dimensional variety, a
linear form $Y_1$ is a primitive element of the
ring extension $\cfq \hookrightarrow \cfq[W]$ if and only if
it separates the points of $W$, i.e. $Y_1(P)\neq
Y_1(Q)$ whenever $P$ and $Q$ are distinct points of $W$.

Let us here remark that this notion of ``geometric solution" has a
long history, which goes back at least to L. Kronecker
\cite{Kronecker82} (see also \cite{Macaulay16}, \cite{Zariski95}).
One might consider \cite{ChGr83} and \cite{GiMo89} as early
references where this notion was implicitly used for the first
time in modern symbolic computation.
%
%
\subsection{Lifting points and lifting fibers}
\label{subsection:lifting_points} Let us consider as in the
previous section an $m$--dimensional $\K$--variety $W$ and a
Noether normalization $\pi:W\to\A^{m}$. We call a point
$P:=(p_1,\ldots,p_{m})\in\A^{m}$ a {\em lifting point} of $\pi$ if
$\pi$ is unramified at $P$, i.e. if the equations
$F_1=0,\ldots,F_{n-m}=0,Y_1=p_1,\ldots,Y_{m}=p_{m}$ define the
fiber $\pi^{-1}(P)$ by transversal cuts. We call the
zero--dimensional variety $W_P:=\pi^{-1}(P)$ the {\em lifting
fiber} of the point $P$.

Suppose that there is given a geometric solution of $W$ as in the
previous section,
and a lifting point $P$ of $\pi$ not vanishing the discriminant of
the polynomial $q^{(m)}$ with respect to the variable $T$. Then
the given geometric solution of the variety $W$ induces a
geometric solution of the lifting fiber $W_P$. This geometric
solution of $W_P$ is given by the linear forms $Y_{m+1}\klk Y_n$,
the polynomial $q^{(m)}(P,T)$ and the parametrization $(\partial
q^{(m)}\!/\partial T)(P,T)Y_{\!m+2}-v^{(m)}_{m+2}(P\!,T), \ldots,
(\partial q^{(m)}\!/\partial T)(P\!,T)Y_{\!n}-
v^{(m)}_{n}\!(P\!,T)$. We call such a geometric solution of $W$
{\em compatible} with the lifting point $P$. \spar

Let us observe that $\pi$ is unramified at a given point
$P\in\A^{m}$ if and only if $J(x)\not=0$ holds for any
$x\in\pi^{-1}(P)$, where $J\in\cfq[\xon]$ denotes the Jacobian
determinant of $\yo{m},\fo{m}$ with respect to the variables
$\xo{n}$. Furthermore, \cite[Proposici\'on 28]{Morais97} shows
$\pi$ is unramified at $P\in\A^{m}$ if and only if the condition
$\#\pi^{-1}(P)=\deg W$ holds.

For $1\le j\le m-n$, let $F_j(Y_1,\ldots,Y_n)$ denote the element
of $\cfq[\yon]$ obtained by rewriting $F_j(\xon)$ in the variables
$\yon$. The following result, probably well--known, is included
here for lack of a suitable reference:
\begin{lemma}
\label{lemma:equiv_ramif} Let notations and assumptions be as
above. Suppose that $\pi$ is unramified at a point $P\in\A^{m}$.
Then the Jacobian matrix $(\partial F_j/\partial Y_{n-k+1})_{1\le
j,k\le n-m}(x)$ is nonsingular for any point $x\in\pi^{-1}(P)$.
\end{lemma}
\begin{proof}
Let $W_P:=\pi^{-1}(P)$, let $\widetilde{\pi}:W_P\to\A^{n-m}$ be
the projection morphism defined by the linear forms $Y_{m+1}\klk
Y_n$ and let $\widetilde{\pi}^*:\cfq[Y_{m+1}\klk Y_n]\to\cfq[W_P]$
denote the corresponding morphism of coordinate rings. Let $I_P$
denote the ideal of $\cfq[Y_{m+1},\ldots,Y_n]$ generated by the
polynomials $F_j(P,Y_{m+1},\ldots,Y_n)$ for $1\le j\le m-n$. We
claim that $I_P$ equals the kernel of the morphism
$\widetilde{\pi}^*$. Indeed, it is clear that the ideal $I_P $ is
included in the kernel of the morphism $\widetilde{\pi}^*$. On the
other hand, let $F\in\cfq[Y_{m+1},\ldots,Y_n]$ satisfy the
condition $\widetilde{\pi}^*(F)=0$. This implies that $F$,
considered as an element of $\cfq[\xon]$, vanishes on any point of
the fiber $W_P$. This implies that the following relation holds:
\begin{equation}\label{eq:ramif_cond} F\in\big(Y_1-p_1\klk
Y_{m}-p_{m},F_1(\yon)\klk F_{n-m}(\yon)\big).\end{equation}
Specializing the variables $Y_1\klk Y_{m}$ into the values
$p_1\klk p_{m}$ in (\ref{eq:ramif_cond}) we conclude that $F\in
I_P$ holds.

From the claim and the fact that $\widetilde{\pi}^*$ is surjective
we deduce the existence of an isomorphism of $\cfq$--algebras:
$$\cfq[\yon]/\big(F_1(P,Y_{m+1},\ldots,Y_n),\ldots,F_{n-m}(P,
Y_{m+1}, \ldots,Y_n) \big)\cong\cfq[W_P].$$
This shows that the ideal $I_P$ is radical. Since $W_P$ is a
zero--dimensional variety, from e.g. \cite[Chapter 4, Corollary
2.6]{CoLiOS98} it follows that $W_P$ is a smooth variety.
Therefore, applying the Jacobian criterion finishes the proof of
the lemma.
\end{proof}
%
%
\subsection{On the algorithmic model}
\label{subsection:slps}
Algorithms in elimination theory are usually described using the
standard dense (or sparse) complexity model, i.e. encoding
multivariate polynomials  by  means  of  the vector of all (or of all nonzero)
coefficients. Taking into account that a generic $n$--variate
polynomial of degree~$d$ has~${d+n\choose n}=O(d^n)$ nonzero
coefficients, we see that the dense or sparse representation of
multivariate polynomials requires an exponential size, and their
manipulation usually requires an exponential number of arithmetic
operations with respect to the parameters~$d$ and~$n$. In order to
avoid this exponential behavior, we are going to use an
alternative encoding  of input, output and intermediate results of
our computations by means of \slps\ (cf. \cite{Heintz89},
\cite{Strassen90}, \cite{Pardo95}, \cite{BuClSh97}). A {\em
\slp}~$\beta$ in~${\K}(\xon)$ is a finite sequence of rational
functions~${(F_1\klk F_k)}\in {\K}(\xon)^k$ such that for~${1\le
i\le k}$, $F_i$ is either an element of the set~$\{\xon\}$, or an
element of~$\K$ (a {\em parameter}), or there exist~${1\le
i_1,i_2<i}$ such that $F_i=F_{i_1}\circ_i\,F_{i_2}$ holds, where
$\circ_i$ is one of the arithmetic operations~${+,-,\times
,{\div}}$. The \slp~$\beta$ is called {\em division--free}
if~$\circ_i$ is different from~${\div}$ for~${1\le i\le k}$. Two
basic natural measures of the complexity of $\beta$ are its {\em
space} and {\em time} (cf. \cite{Borodin93}, \cite{Savage98}).
Space is defined as the maximum number of arithmetic registers
used in the evaluation process defined by $\beta$, and time is
defined as the total number of arithmetic operations performed
during the evaluation. We say that the \slp~$\beta$ {\em computes}
or {\em represents} a subset~$S$ of~${\K}(\xon)$ if~$S\subset
\{F_1\klk F_k\}$ holds.

Our model of computation is based on the concept of \slps.
However, a model of computation consisting {\em only} of \slps\ is
not expressive enough for our purposes. Therefore we allow our
model to include decisions and selections (subject to previous
decisions). For this reason we shall also consider {\em
computation trees}, which are \slps\ with {\em branchings}. Time
and space of the evaluation of a given computation tree are
defined analogously as in the case of \slps\ (see e.g.
\cite{Gathen86}, \cite{BuClSh97} for more details on the notion of
computation trees). \spar

A difficult point in the manipulation of multivariate polynomials
over finite fields is the so--called {\em identity testing
problem}: given two elements $F$ and $G$ of ${\K}[\xon]$, decide
whether $F$ and $G$ represent the same polynomial function on
${\K}^n$. Indeed, all known deterministic algorithms solving this problem have
complexity at least~$(\#\K)^{\Omega(1)}$. In this article we are
going to use {\em probabilistic} algorithms to solve the identity
testing problem, based on the following result:
\begin{theorem}[\cite{LiNi83}, \cite{Schmidt76}]
\label{th:Zippel_Schwartz} Let $F$ be a nonzero polynomial
of~$\cfq[\xon]$ of degree at most~$d$ and let $\K$ be a finite
field extension of $\fq$. Then the number of zeros of $F$ in
${\K}^n$ is at most $d(\# {\K})^{n-1}$.
\end{theorem}
For the  analysis of our algorithms, we  shall interpret the
statement of Theorem~\ref{th:Zippel_Schwartz} in terms  of
probabilities. More precisely, given a fix nonzero polynomial $F$
in $\cfq[\xon]$ of degree at most~$d$, from
Theorem~\ref{th:Zippel_Schwartz} we conclude that the probability
of choosing  randomly a  point~$a\in {\K}^n$ such that~${F(a)=0}$
holds is bounded from above by~$d/ \#{\K}$ (assuming a
uniform distribution of probability on the elements of ${\K}^n$).
%
%
\section{On the preparation of the input data}
\label{section:preparation}

From now on, let $n\ge 3$ and $d\ge 2$, and let be given
polynomials $F_1\klk F_r\in \fq[X_1,\dots,\!X_n]$ of maximum
degree $d$, which generate a radical ideal of $\fq[\xon]$ and form
a regular sequence of $\fq[\xon]$. Suppose further that $\fo s$ generate a
radical ideal for $1\le s\le r-1$ and $V_r:=V(\fo{r})$ is absolutely
irreducible.

In the sequel we shall consider algorithms which ``solve"
symbolically the (input) equation system $F_1=0,\ldots,F_r=0$ over
$\cfq$. As in \cite{GiHeMoMoPa98} and \cite{GiHaHeMoMoPa97}, we
associate to the equation system $F_1=0,\ldots,F_r=0$ a parameter
$\delta$, called the (geometric) {\em degree} of the system, which
is defined as follows: for $1\le s\le r$, let $V_s\subset\A^n$ be
the $\fq$--variety defined by $\fo s$ and let $\delta_s$ denote
its degree. The geometric degree of the system
$F_1=0,\ldots,F_r=0$ is then defined as $\delta:=\max_{1\le s\le
r}\delta_s.$

In this section we are going to determine a genericity condition
underlying the choice of a simultaneous Noether normalization of
the varieties $V_1\klk V_r$ and lifting points $P^{(s)}\in
\A^{n-s}$ ($1\le s\le r$) such that, for $1\le s\le r-1$, the
lifting fiber $V_{P^{(s+1)}}$ has the following property: for any
point $P\in V_{P^{(s+1)}}$, the morphism $\pi_s$ is unramified at
$\pi_s(P)$. By a simultaneous Noether normalization we understand
a linear change of variables such that the new variables
$Y_1,\ldots,Y_n$ are in Noether position with respect to $V_s$ for
$1\le s\le r$. Finally, we are going to find an
affine linear subspace $L$ of $\A^n$ of dimension $r+1$ such that
$V_r\cap L$ is an absolutely irreducible curve of $\A^n$ of degree
$\delta_r$.
%
%
\subsection{Simultaneous Noether normalization}
\label{subsec:Noether} It is well--known that a generic choice of
linear forms $Y_1\klk Y_n$ yields a simultaneous Noether
normalization of the varieties $V_1\klk V_r$. In order to prove
the existence of a simultaneous Noether normalization defined over
a given finite field extension of $\fq$ we need suitable
genericity conditions. The following proposition yields an upper bound
on the degree of the genericity condition underlying the choice of
such linear forms.

\begin{proposition}
\label{prop:chow} Let us fix $s$ with $1\le s\le r$. Let
$\Lambda:=(\Lambda_{ij})_{1\le i\le n-s+1,1\le j\le n}$ be a
matrix of indeterminates, let
$\Lambda^{(i)}:=(\Lambda_{i1},\dots,\Lambda_{in})$ for $1\le i\le
n-s+1$ and let $\Gamma:=(\Gamma_1\klk \Gamma_{n-s+1})$ be a vector
of indeterminates. Let $X:=(\xon)$ and let $\widetilde{Y}:=\Lambda
X+\Gamma$. Then there exists a nonzero polynomial
$A_s\in\cfq[\Lambda,\Gamma]$ of degree at most
$2(n-s+2)\delta_s^2$ such that, for any
$(\lambda,\gamma)\in\A^{(n-s+1)n}\times \A^{n-s+1}$ with
$A_s(\lambda,\gamma)\not=0$, the following conditions are
satisfied:
\begin{itemize}
  \item[$(i)$]
Let $Y:=\lambda X+\gamma:=
  (Y_1\klk Y_{n-s+1})$. Then the mapping $\pi_s:V_s\to \A^{n-s}$ defined
  by $Y_1\klk Y_{n-s}$ is a finite morphism.
  \item[$(ii)$] The linear form $Y_{n-s+1}$ induces a primitive element
  of the integral ring extension $R_s:=\cfq[\yo{n-s}]\hookrightarrow
  \cfq[V_s]$.
\end{itemize}
\end{proposition}

\begin{proof}
Let us consider the following morphism of algebraic varieties:
\begin{equation}\label{eq:morph_chow}\begin{array}{crcl}
  \Phi:& \A^{(n-s+1)n}\times\A^{n-s+1}\times V_s & \to &
  \A^{(n-s+1)n}\times\A^{n-s+1}\times \A^{n-s+1} \\
       & (\lambda,\gamma,x) & \mapsto & (\lambda,\gamma,\lambda x+\gamma)
\end{array}\end{equation}
Using standard facts about Chow forms (see e.g.
\cite{Shafarevich84}, \cite{Caniglia90}, \cite{KrPaSo01}), we
deduce that $\overline{Im(\Phi)}$ is a hypersurface of
$\A^{(n-s+1)n}\times\A^{n-s+1}\times \A^{n-s+1}$, defined by a
squarefree polynomial
$P_{V_s}\in\cfq[\Lambda,\Gamma,\widetilde{Y}_1\klk
\widetilde{Y}_{n-s+1} ]$ which satisfies the following degree
estimates:
\begin{itemize}
  \item $\deg_{\widetilde{Y}}\!P_{V_s}=
  \deg_{\widetilde{Y}_{n-s+1}}\!P_{V_s}= \delta_s$,
  \item $\deg_{\Lambda^{(i)}\!,\,\Gamma_i}\!P_{V_s}\le
  \delta_s$ for $1\le i\le n-s+1$.
\end{itemize}

Let $A_{1,s}\in\cfq[\Lambda,\Gamma]$ be the (nonzero) polynomial
which arises  as coefficient of the monomial
$\widetilde{Y}_{n-s+1}^{\delta_s}$ in the polynomial $P_{V_s}$,
considering $P_{V_s}$ as an element of
$\cfq[\Lambda,\Gamma][\widetilde{Y}]$. The above estimates imply
$\deg A_{1,s}\le (n-s+1)\delta_s$. Let
$\widetilde{A}_{1,s}\in\cfq[\Lambda^{\!(i)},\Gamma_i:1\le i\le
n-s]$ be a nonzero polynomial arising as coefficient of a monomial
of $A_{1,s}$, considering $A_{1,s}$ as an element of
$\cfq[\Lambda^{(i)},\Gamma_i:1\le i\le
n-s][\Lambda^{(n-s+1)},\Gamma_{n-s+1}]$.

Let $(\lambda^*,\gamma^*)\in\A^{(n-s)n}\times \A^{n-s}$ be any
point for which $\widetilde{A}_{1,s}(\lambda^*,\gamma^*)\not=0$
holds, and let $Y:=(Y_1\klk Y_{n-s}):=\lambda^* X+\gamma^*$. We
claim that condition $(i)$ of the statement of Proposition
\ref{prop:chow} holds. Indeed, since
$A_{1,s}^*:=A_{1,s}(\lambda^*,\gamma^*,\Lambda^{\!(n-s+1)},\Gamma_{n-s+1})$
is a nonzero element of $\cfq[\Lambda^{\!(n-s+1)},\Gamma_{n-s+1}]
$, we deduce the existence of $n$ $\cfq$--linear independent
vectors $w_1\klk w_n\in\A^n$ and values $a_1\klk a_n\in\A^1$ such
that $A_{1,s}^*(w_j,a_j)\not=0$ holds for $1\le j\le n$. Let
$\ell_j:=w_j X+a_j$ for $1\le j\le n$. By construction, for $1\le
j\le n$ the polynomial
$P_{V_s}(\lambda^*,\gamma^*,w_j,a_j,\yo{n-s},\ell_j)$ is an
integral dependence equation for the coordinate function induced
by  $\ell_j$ in the ring extension $R_s\hookrightarrow\cfq[V_s]$.
Since $\cfq[\ell_1\klk\ell_n]=\cfq[\xon]$, we conclude that
condition $(i)$ holds. \spar

Furthermore, since $\cfq[\Lambda,\Gamma,\widetilde{Y}]/(P_{V_s})$
is a reduced $\cfq$--algebra and $\cfq$ is a perfect field, from
\cite[Proposition 27.G]{Matsumura80} we conclude that the
(zero--dimensional)
$\cfq(\Lambda,\Gamma,\widetilde{Y}_1\klk\widetilde{Y}_{\!n-s})$--algebra
$\cfq(\Lambda,\Gamma,\widetilde{Y}_1\klk\widetilde{Y}_{\!n-s})
[\widetilde{Y}_{\!n-s+1}]/(P_{V_s})$ is reduced. This implies that
$P_{V_s}$ is a separable element of
$\cfq(\Lambda,\Gamma,\widetilde{Y}_1\klk
\widetilde{Y}_{n-s})[\widetilde{Y}_{n-s+1}]$ and hence $P_{V_s}$
and $\partial P_{V_s}/\partial \widetilde{Y}_{n-s+1}$ are
relatively prime in
$\cfq(\Lambda,\Gamma,\widetilde{Y}_1\klk\widetilde{Y}_{n-s})
[\widetilde{Y}_{n-s+1}]$. Then the discriminant
\begin{equation}\label{eq:def_disc_chow}
\rho_s:=\mbox{Res}_{\widetilde{Y}_{n-s+1}}(P_{V_s},\partial
P_{V_s} /\partial \widetilde{Y}_{n-s+1})\end{equation}
of $P_{V_s}$ with respect to $\widetilde{Y}_{n-s+1}$ is a nonzero
element of $\cfq[\Lambda,\Gamma,\widetilde{Y}_1\klk
\widetilde{Y}_{n-s}]$ which satisfies the following degree
estimates:
\begin{itemize}
  \item $\deg_{\widetilde{Y}_1\!\klk\,\widetilde{Y}_{n-s}}\rho_s\le
(2\delta_s-1)\delta_s$.
  \item $\deg_{\Lambda^{(i)}\!,\,\Gamma_i}\rho_s\le
(2\delta_s-1)\delta_s$ for $1\le i\le n-s+1$.
\end{itemize}

Let $\rho_{1,s}\in \cfq[\Lambda,\Gamma]$ a nonzero coefficient of
a monomial of $\rho_s$, considering $\rho_s$ as an element of
$\cfq[\Lambda,\Gamma][\widetilde{Y}_1\klk \widetilde{Y}_{n-s}]$,
and let $A_s:=\rho_{1,s}\widetilde{A}_{1,s}$. Observe that $\deg
A_s\le 2(n-s+2)\delta_s^2$ holds. Let
$(\lambda,\gamma)\in\A^{(n-s+1)n} \times\A^{n-s+1}$ satisfy the
condition $A_s(\lambda,\gamma)\not=0$, let $Y:=\lambda X+\gamma$
and denote by $(\lambda^*,\gamma^*)\in\A^{(n-s)n} \times\A^{n-s}$
be the matrix formed by the first $n-s$ rows of
$(\lambda,\gamma)$. Let $\rho_s^*$ be the polynomial obtained from
$\rho_s$ by specializing the variables $\Lambda^{(i)},\Gamma_i$
($1\le i\le n-s$) into the value $(\lambda^*,\gamma^*)$. Then
$\rho_s^*$ is a nonzero element of
$\cfq[\Lambda^{(n-s+1)},\Gamma_{n-s+1},\yo{n-s}]$ which equals the
discriminant of $P_{V_s}(\lambda^*,\Lambda^{\!(n-s+1)}, \gamma^*,
\Gamma_{n-s+1},\yo{n-s},\widetilde{Y}_{n-s+1})$ with respect to
$\widetilde{Y}_{n-s+1}$. It is clear that condition $(i)$ holds.
We claim that condition $(ii)$ holds. \spar

Let $\xi_1\klk \xi_n$ be the coordinate functions of $V_s$ induced
by $\xon$, let $\zeta_i:=\sum_{j=1}^n\lambda_{i,j}\xi_j+\gamma_i$
for $1\le i\le n-s$ and let $\widehat{Y}_{n-s+1}:=\sum_{j=1}^n
\Lambda_{n-s+1,j}\xi_j+\Gamma_{n-s+1}$. From the properties of the
Chow form of $V_s$ we conclude that the identity
\begin{equation}
\label{equation:chow1}
\begin{array}{ccl}
0\!\!\!\!\!\!&=&\!\!\!\!\!\!P_{V_s}\!(\lambda^*\!,\!\Lambda^{\!(n-s+1)}\!,
\!\gamma^*\!,
\Gamma_{\!n-s+1},\!\zeta_1\klk\!\zeta_{n-s},\widehat{Y}_{n-s+1})\\
\!\!&=&\!\!\!\!\!\!P_{V_s}\!(\lambda^*\!,\!\Lambda^{\!(n-s+1)}\!,\!\gamma^*\!,
\Gamma_{\!n-s+1},\!\zeta_1,\ldots,\!\zeta_{n-s},\Lambda_{n-s+1,1}\xi_1\!\plp
\!\Lambda_{n-s+1,n}\xi_n)
\end{array}
\end{equation}
holds in $\cfq[\Lambda^{\!(n-s+1)},\Gamma_{n-s+1}]\otimes_{\cfq}
\cfq[V_s]$.
Following e.g. \cite{AlBeRoWo96} or \cite{Rouillier97}, from
(\ref{equation:chow1}) we deduce that the following identity holds
in $\cfq[\Lambda^{(n-s+1)},\Gamma_{n-s+1}]\otimes_{\cfq}
\cfq[V_s]$ for $1\le k\le n$:
\begin{equation}
\label{equation:chow2}
\begin{array}{l}
(\partial P_{V_s}/\partial \widetilde{Y}_{n-s+1})
(\lambda^*,\Lambda^{(n-s+1)},\gamma^*,\Gamma_{n-s+1},\zeta_1 \klk
\zeta_{n-s},\widehat{Y}_{n-s+1})\xi_k+ \\\  +(\partial
P_{V_s}/\partial \Lambda_{n-s+1,k})
(\lambda^*,\Lambda^{(n-s+1)},\gamma^*,\Gamma_{n-s+1},\zeta_1\klk
\zeta_{n-s},\widehat{Y}_{n-s+1})=0.
\end{array}
\end{equation}
Since $\rho_s^*(\Lambda^{\!(n-s+1)},\Gamma_{n-s+1},Y_1\klk
Y_{n-s})$ is the discriminant of the polynomial\linebreak
$P_{V_s}(\lambda^*\!,\Lambda^{\!(n-s+1)},\gamma^*\!,\Gamma_{n-s+1},\yo{n-s},
\widetilde{Y}_{n-s+1})$ with respect to $\widetilde{Y}_{n-s+1}$,
it can be written as a linear combination of
$P_{V_s}\!(\lambda^*\!,\Lambda^{\!(n-s+\!1)}\!,\!\gamma^*\!,\!\Gamma_{\!n-s+\!1},
\!\yo{\!n-s},\!\widetilde{Y}_{n-s+1}\!)$ and $(\partial
P_{V_s}/\partial \widetilde{Y}_{n-s+1}) (\lambda^*,
\Lambda^{(n-s+1)},
\gamma^*,\Gamma_{n-s+1},\yo{n-s},\widetilde{Y}_{n-s+1})$.
Combining this observation with (\ref{equation:chow1}) and
(\ref{equation:chow2}) we conclude that
\begin{equation}
\label{equation:chow3}
\begin{array}{c}
\rho_s^*(\Lambda^{\!(n-s+1)},\Gamma_{n-s+1},\zeta_1\klk
\zeta_{n-s})\xi_k+
\qquad\qquad\qquad\qquad\qquad\qquad\\\qquad\qquad\qquad
\qquad+P_k (\Lambda^{\!(n-s+1)},\Gamma_{n-s+1},\zeta_1\klk
\zeta_{n-s},\widehat{Y}_{n-s+1})=0
\end{array}\end{equation}
holds, where $P_k$ is a nonzero element of
$\cfq[\Lambda^{\!(n-s+1)},\Gamma_{n-s+1},\om Z {n-s+1}]$ for $1\le
k\le n$. Specializing identity (\ref{equation:chow3}) into the
values $\Lambda_{n-s+1,j}:=\lambda_{n-s+1,j}$ ($1\le j\le n$) and
$\Gamma_{n-s+1}=\gamma_{n-s+1}$ for $1\le k\le n$, we conclude
that $Y_{n-s+1}$ induces a primitive element of the
$\cfq$--algebra extension $\cfq(Y_1\klk Y_{n-s})\hookrightarrow
\cfq(Y_1\klk Y_{n-s})\otimes_{\cfq}\cfq[V_s]$.

Condition ($i$) implies that $\cfq[V_s]$ is a finite free
$R_s:=\cfq[Y_1\klk Y_{n-s}]$--module and hence $\cfq(Y_1\klk
Y_{n-s})\otimes_{\cfq}\cfq[V_s]$ is a finite--dimensional
$\cfq(Y_1\klk Y_{n-s})$--vector space. Furthermore, the dimension
of $\cfq(Y_1,\dots,\!Y_{n-s})\otimes_{\cfq}\cfq[V_s]$ as
$\cfq(Y_1\klk Y_{n-s})$--vector space equals the rank of
$\cfq[V_s]$ as $R_s$--module. On the other hand, since $R_s$ is
integrally closed, we have that the minimal dependence equation of
an arbirtrary element $f\in\cfq[V_s]$ over $\cfq(Y_1\klk Y_{n-s})$ equals
the minimal integral dependence equation of $f$ over $R_s$ (see e.g.
\cite[Lemma II.2.15]{Kunz85}). Combining this remark with the fact
that $Y_{n-s+1}$ induces a primitive element of the
$\cfq$--algebra extension $\cfq(Y_1\klk
Y_{n-s})\hookrightarrow\cfq(Y_1\klk Y_{n-s}) \otimes_{\cfq}
\cfq[V_s]$, we conclude that $Y_{n-s+1}$ also induces a primitive
element of the $\cfq$--algebra extension
$R_s\hookrightarrow\cfq[V_s]$. This shows condition ($ii$) and
finishes the proof of the proposition.
\end{proof}
%
%
\subsection{Lifting fibers not meeting a discriminant}
\label{subsec:main_preparation} Our second step is to find lifting
points $P^{(s+1)}\in\A^{n-s-1}$ for $0\le s \le r-1$ such that
the corresponding lifting fiber $V_{P^{(s+1)}}$ has the following
property: for any point $P\in V_{P^{(s+1)}}$, the morphism $\pi_s$
is unramified at $\pi_s(P)$. With this condition we shall be able
to find a geometric solution of the variety $V_s$ such that no point
$P\in V_{P^{(s+1)}}$ annihilates the discriminant of the corresponding
minimal polynomial $q^{(s)}$, which in turn will allow us to
avoid dealing with multiplicities during the computations.

For this purpose we need the following technical result, which is
a slightly simplified version of \cite[Lemma 1 $(iii)$]{HeMaWa01}
with an improved degree estimate.
\begin{lemma}\label{lemma:degree_estimate} With notations and
assumptions as above, let us fix $s$ with $1\le s\le r$. Let $A_s$
be the polynomial of the statement of Proposition \ref{prop:chow}
and let be given a polynomial $H\in\cfq[\Lambda,\Gamma,X]$ of
degree at most $D$. Suppose that the Zariski closure
$\widehat{V}_s$ of the set $(\A^{(n-s+1)n}\times\A^{n-s+1}\times
V_s)\cap\{H=0,A_s\not=0\}$ satisfies the condition
$\dim\,\widehat{V}_s\le (n-s+1)(n+2)-2$. Then the Zariski closure
of the image of $\widehat{V}_s$ under the morphism
$\Phi^*:\A^{(n-s+1)n}\times\A^{n-s+1}\times
V_s\to\A^{(n-s+1)n}\times\A^{n-s+1}\times \A^{n-s}$ defined by
$\Phi^*(\lambda,\gamma,x):=(\lambda,\gamma,\lambda^*x+\gamma^*)$
is empty or is contained in a hypersurface of
$\A^{(n-s+1)n}\times\A^{n-s+1}\times\A^{n-s}$ of degree at most
$2(n-s+2)D\delta_s^2$ (here $\lambda^*$ and $\gamma^*$ denote the
first $n-s$ rows of $\lambda$ and $\gamma$ respectively).
\end{lemma}
\begin{proof}
We use the notations of the proof of Proposition \ref{prop:chow}.
Since the Chow form $P_{V_s}$ of the variety $V_s$ is a separable
element of $\cfq(\Lambda,\Gamma,\widetilde{Y}_1,\dots,\!
\widetilde{Y}_{n-s})[\widetilde{Y}_{n-s+1}]$, we conclude that
$\partial P_{V_s}/\partial \widetilde{Y}_{n-s+1}$ is not a zero
divisor of $\cfq[\Lambda,\Gamma,\widetilde{Y}]/(P_{V_s})$, and
hence of the $\cfq$--algebra $\cfq[\Lambda,\Gamma]\otimes_{\cfq}
\cfq[V_s]$. Arguing as in the proof of identity
(\ref{equation:chow2}), we see that the following identity holds
in $\cfq[\Lambda,\Gamma] \otimes_{\cfq} \cfq[V_s]$:
\begin{equation}
\label{equation:chow4} (\partial P_{V_s}/\partial
\widetilde{Y}_{n-s+1})
(\Lambda,\Gamma,\widehat{Y})\,\xi_k+(\partial P_{V_s}/\partial
\Lambda_{n-s+1,k}) (\Lambda,\Gamma,\widehat{Y})=0,
\end{equation}
where $\widehat{Y}:=\Lambda\xi+\Gamma$ and $\xi:=(\xi_1\klk\xi_n)$
is the vector of coordinate functions of $V_s$ induced by $X$.

Let $\widehat{H}\in\cfq[\Lambda,\Gamma,\widetilde{Y}]$ be the
polynomial obtained by replacing in $H$ the variable $X_k$ by
$-(\partial P_{V_s}/\partial \widetilde{Y}_{n-s+1})^{-1}(\partial
P_{V_s}/\partial \Lambda_{n-s+1,k})$ for $1\le k\le n$ and
clearing denominators. Observe that
$\deg_{\widetilde{Y}}\widehat{H}=
\deg_{\widetilde{Y}_{n-s+1}}\widehat{H}\le D\delta_s$ and
$\deg_{\Lambda,\Gamma}\widehat{H}\le (n-s+1)D\delta_s$ holds.

Let $R:=\mbox{Res}_{\widetilde{Y}_{n-s+1}}
(P_{V_s},\widehat{H})\in\cfq[\Lambda,
\Gamma,\widetilde{Y}_1\klk\widetilde{Y}_{n-s}]$ be the resultant
of $P_{V_s}$ and $\widehat{H}$ with respect to the variable
$\widetilde{Y}_{n-s+1}$. Observe that the Sylvester matrix of
$P_{V_s}$ and $\widehat{H}$ is a matrix of size at most
$(D+1)\delta_s\times(D+1)\delta_s$ with at most $D\delta_s$
columns consisting of coefficients of $P_{V_s}$ or zero entries,
and at most $\delta_s$ columns consisting of coefficients of
$\widehat{H}$ or zero entries. This shows that $\deg R\le
2(n-s+2)D\delta_s^2$ holds. On the other hand, from identity
(\ref{equation:chow4}) and the properties of the resultant we
conclude that
$R(\Lambda,\Gamma,\widetilde{Y}_1\klk\widetilde{Y}_{n-s})$
vanishes on the variety $\widehat{V}_s$. Furthermore, the
assumption $\dim \,\widehat{V}_s\le (n-s+1)(n+2)-2$ implies
$R(\Lambda,\Gamma,\widetilde{Y}_1\klk\widetilde{Y}_{n-s}) \not=0$.
This finishes the proof of the lemma.
\end{proof}

Now we are ready to prove the main theorem of this section. This
result states an appropriate bound for the degree of a certain
polynomial, whose nonvanishing expresses a suitable genericity
condition for the coefficients of the linear forms \yon\ and the
coordinates of the lifting points $P^{(s+1)}$ ($1\le s\le r-1$) we
are looking for. Let us remark that a similar result, with higher
degree estimates, is proved in \cite[Theorem 3]{HeMaWa01} for a
$\Q$--definable affine equidimensional variety of $\C^n$.
Unfortunately, the proof of \cite[Theorem 3]{HeMaWa01} makes an
essential use of the fact that the underlying variety is defined
over $\Q$ and therefore cannot be used in our situation. On the
other hand, we obtain a significant improvement of the degree
estimates of \cite[Theorem 3]{HeMaWa01}, which is a critical point
for our subsequent purposes.

\begin{theorem}
\label{theorem:simNoether} Let notations be as in Proposition
\ref{prop:chow} and let us fix $s$ with $1\le s<r$. Then there
exists a nonzero polynomial $B_s\in\cfq[\Lambda,
\Gamma,\widetilde{Y}_1\klk\widetilde{Y}_{n-s}]$, of degree at most
$4(n-s+3)^2nd\delta_s^2\delta_{s+1}^2$, such that for any
$(\lambda,\gamma,P)\in\A^{(n-s+1)n}\times \A^{n-s+1}\times
\A^{n-s}$ with $B_s(\lambda,\gamma,P)\not=0$ the following
conditions are satisfied:
\begin{itemize}
\item[$(i)$] Let $Y:=(\yo{n-s+1}):=\lambda X+\gamma$. Then the
mapping $\pi_s:V_s\to\A^{n-s}$ defined by
$\pi_s(x):=\big(Y_1(x)\klk Y_{n-s}(x)\big)$ is a finite morphism,
$P\in\A^{n-s}$ is a lifting point of $\pi_s$ and $Y_{n-s+1}$ is a
primitive element of $\pi_s^{-1}(P)$.
\item[$(ii)$] Let $P^*\in\A^{n-s-1}$ be the vector that consists of the
first $n-s-1$ coordinates of $P$. Then the mapping
$\pi_{s+1}:V_{s+1} \to \A^{n-s-1}$ defined by
$\pi_{s+1}(x):=\big(Y_1(x)\klk Y_{n-s-1}(x)\big)$ is a finite
morphism, $P^*$ is a lifting point of $\pi_{s+1}$ and $Y_{n-s}$ is
a primitive element of $\pi_{s+1}^{-1}(P^*)$.
\item[$(iii)$] Any point $Q\in\pi_s\big(\pi_{s+1}^{-1}(P^*)\big)$ is
a lifting point of $\pi_s$ and $Y_{n-s+1}$ is a primitive element
of $\pi_s^{-1}(Q)$ for any $Q\in\pi_s\big(\pi_{s+1}^{-1}(P^*)\big)$.
\end{itemize}
\end{theorem}

\begin{proof}
Let $A_s$ and $A_{s+1}$ be the polynomials obtained by applying
Proposition \ref{prop:chow} to the varieties $V_s$ and $V_{s+1}$
respectively. Let $D_s,D_{s+1}\in\cfq[\Lambda,\Gamma,X]$ be the
following polynomials: $$D_s:=\det\begin{pmatrix}
  \Lambda_{1,1} & \dots & \Lambda_{1,n} \\
  \vdots &  & \vdots \\
  \Lambda_{n-s,1} & \dots & \Lambda_{n-s,n} \\
  \frac{\partial F_1}{\partial X_1} & \dots &
  \frac{\partial F_1}{\partial X_n} \\
  \vdots &  & \vdots \\
    \frac{\partial F_s}{\partial X_1} & \dots &
    \frac{\partial F_s}{\partial X_n}
\end{pmatrix},\quad
D_{s+1}:=\det\begin{pmatrix}
  \Lambda_{1,1} & \dots & \Lambda_{1,n} \\
  \vdots &  & \vdots \\
  \Lambda_{n-s-1,1} & \dots & \Lambda_{n-s-1,n} \\
  \frac{\partial F_1}{\partial X_1} & \dots &
  \frac{\partial F_1}{\partial X_n} \\
  \vdots &  & \vdots \\
    \frac{\partial F_{s+1}}{\partial X_1} & \dots &
    \frac{\partial F_{s+1}}{\partial X_n}
\end{pmatrix}.$$

We claim that the Zariski closure of the set $(\A^{(n-s+1)n}
\times \A^{n-s+1}\times V_s)\cap\{D_s=0,A_s\not=0\}$ is empty or
an equidimensional affine subvariety of $\A^{(n-s+1)n}\times
\A^{n-s+1}\times\A^n$ of dimension $(n-s+1)(n+2)-2$. \spar

In order to prove this claim, let
$V_s=\mathcal{C}_1\cup\cdots\cup\mathcal{C}_t$ be the
decomposition of $V_s$ into irreducible components. Then we have
that $\A^{(n-s+1)n}\times\A^{n-s+1}\times V_s=\cup_{i=1}^t
\A^{(n-s+1)n}\times\A^{n-s+1}\times\mathcal{C}_i$ is the
decomposition of $\A^{(n-s+1)n}\times\A^{n-s+1}\times V_s$ into
irreducible components. Let $\A^{(n-s+1)n}\times\A^{n-s+1}
\times\mathcal{C}$ be any of these irreducible components and let
$x\in\mathcal{C}$ be a nonsingular point of $V_s$. Then
$D_s(\Lambda,x)\not=0$ holds and therefore there exists
$\lambda\in\A^{(n-s+1)n}$ such that $D_s(\lambda,x)\not=0$ holds.
This shows that there exists a point
$(\lambda,\gamma,x)\in\A^{(n-s+1)n}\times\A^{n-s+1}\times
\mathcal{C}$ not belonging to the hypersurface $\{D_s=0\}$. On the
other hand, $D_s(0,x)=0$ holds for any $x\in V_s$, where $0$
represents the zero matrix of $\A^{(n-s+1)n}$. This proves that
$\{D_s=0\}\cap(\A^{(n-s+1)n}\times\A^{n-s+1}\times V_s)$ is an
equidimensional variety of dimension $(n-s+1)(n+2)-2$ and hence
the Zariski closure of the set $(\A^{(n-s+1)n}\times\A^{n-s+1}
\times V_s) \cap\{D_s=0,A_s\not=0\}$ is either empty or an
equidimensional variety of dimension $(n-s+1)(n+2)-2$. This shows
the claim.\spar

The same argument {\em mutatis mutandis} shows that the Zariski
closure of the set $(\A^{(n-s)n}
\times \A^{n-s}\times V_{s+1})\cap\{D_{s+1}=0,A_{s+1}\not=0\}$ is
empty or an equidimensional affine subvariety of $\A^{(n-s)n}\times
\A^{n-s}\times\A^n$ of dimension $(n-s)(n+2)-2$
\spar

Let us consider the following morphisms:
\[\begin{array}{r}
  \Phi_s: (\A^{(n-s+1)n}\!\times\!
\A^{n-s+1}\!\times\! V_s)\cap\{D_s=0,A_s\not=0\}\to
\A^{(n-s+1)n}\!\times\!\A^{n-s+1}\!\times\!\A^{\!n-s}
\\(\lambda,\gamma,x)\mapsto\big(\lambda,\gamma,
Y_1(x)\klk Y_{n-s}(x)\big)\
\end{array}\]

\[\begin{array}{r}
\Phi_{\!s+\!1}\!:\!(\A^{(n-\!s)n}\!\times\!
\A^{n-\!s}
\!\times\!V_{\!s+\!1})\cap\{D_{\!s+\!1}\!=0,A_{s+\!1}
\!\not=0\}\to
\A^{(n-s)n}\times\A^{n-s}\times\A^{n-s-1}\\
\quad\qquad\qquad\qquad\qquad\qquad\qquad\qquad\qquad
(\lambda^*,\gamma^*,x)\mapsto \!\big(\!\lambda^{*}\!,
\gamma^{*}\!,\!Y_1(x),\dots,\!Y_{\!n-\!s-\!1}(x)\big)
\end{array}\]
From the claims above and Lemma \ref{lemma:degree_estimate} we
deduce that the Zariski closure of $Im(\Phi_s)$ is contained in a
hypersurface of $\A^{(n-s+1)n}\times\A^{n-s+1}\times\A^{n-s}$ of
degree at most $2(n-s+2)n(d-1)\delta_s^2$, and the Zariski closure
of $Im(\Phi_{s+1})$ is contained in a hypersurface of
$\A^{(n-s)n}\times\A^{n-s}\times\A^{n-s-1}$ of degree at most
$2(n-s+1)n(d-1)\delta_{s+1}^2$. Let us denote by
$\widehat{B}_s\in\cfq[\Lambda,\Gamma,\widetilde{Y}_1\klk\widetilde{Y}_{n-s}]$
and $\widehat{B}_{s+1}\in\cfq[\Lambda,\Gamma,
\widetilde{Y}_1\klk\widetilde{Y}_{n-s-1}]$ the polynomials
defining these hypersurfaces respectively.

Let $\rho_s,\rho_{s+1}\in\cfq[\Lambda,\Gamma,\widetilde{Y}_1\klk
\widetilde{Y}_{n-s}]$ be the (nonzero) discriminants of
the varieties $V_s$ and $V_{s+1}$, as defined in eq. (\ref{eq:def_disc_chow})
of the proof of Proposition \ref{prop:chow}.
Recall that
%
$\deg\rho_s\le
(n-s+2)(2\delta_s^2-\delta_s)$ and $\deg\rho_{s+1}\le
(n-s+1)(2\delta_{s+1}^2-\delta_{s+1})$ holds. \spar

\begin{claim} The locally closed set $(\A^{(n-s+1)n}\times \A^{n-s+1}
\times V_{s+1})\cap\{\rho_s\widehat{B}_s=0,A_{s+1}\not=0\}$ has
dimension at most $(n-s+1)(n+2)-3$. \end{claim}

\noindent {\it Proof of Claim.} Let us observe that, from the
definition of the polynomial $A_s$, we deduce that the mapping
$\Phi_s$ above induces a finite morphism of varieties, which we
shall denote also $\Phi_s$ by a slight abuse of notation:
$$\begin{array}{r}
 \!\!\!\Phi_s\!:\!(\A^{\!(n-\!s+\!1)n}\!\times\!\A^{\!n-\!s+\!1}\!
 \times\!V_s)
  \!\cap\!\{A_s\!\not=0\}\!\to\!\big((\A^{(n-\!s+\!1)n}\!\times\!
  \A^{n-\!s+\!1})\!\cap\!\{A_s\!\not=0\}\big)\!\times\!\A^{n-\!s}
\\ (\lambda,\gamma,x)\mapsto\big(\lambda,\gamma,Y_1(x)\klk
Y_{n-s}(x)\big)\qquad\qquad\qquad \end{array}$$
Since $(\A^{(n-s+1)n}\times \A^{n-s+1}\times V_s)\cap\{D_s=0,A_s
\not=0\}$ is an equidimensional subvariety of $(\A^{(n-s+1)n}
\times \A^{n-s+1}\times V_s)\cap\{A_s\not=0\}$ of dimension
$(n-s+2)(n+1)-2$, we see that $\Phi_s(\{D_s=0\})$ is a
hypersurface of $(\A^{(n-s+1)n}\times\A^{n-s+1}
\cap\{A_s\not=0\})\times\A^{n-s}$, which is therefore definable by
the polynomial $\widehat{B}_s$. This means that the identity
$\Phi_s(\{D_s=0,A_s\not=0\})=\{\widehat{B}_s=0, A_s\not=0\}$
holds.

From the cylindrical structure of the variety $\A^{(n-s+1)n}\times
\A^{n-s+1}\times V_{s+1}$ we conclude that no irreducible
component of this variety is contained in $\{A_s=0\}$. This
implies that $\mathcal{D}\cap\{A_s\not=0\}$ is a dense open subset
of $\mathcal{D}$ for any irreducible component $\mathcal{D}$ of
$\A^{(n-s+1)n}\times \A^{n-s+1}\times V_{s+1}$. Suppose that there
exists an irreducible component $\mathcal{D}$ of
$\A^{(n-s+1)n}\times \A^{n-s+1}\times V_{s+1}$ contained in
$\Phi_s^{-1}(\{\rho_s\widehat{B}_s=0\})$. Then
$$\mathcal{D}\cap\{A_s\not=0\}\subset\Phi_s^{-1}(\{\rho_s\widehat{B}_s=0\})
\cap\{A_s\not=0\}= \Phi_s^{-1}(\{\rho_s\widehat{B}_s=0\}
\cap\{A_s\not=0\}),$$
which implies
$$\Phi_s(\mathcal{D}\cap\{A_s\not=0\})\subset \Phi_s\circ
\Phi_s^{-1}(\{\rho_s\widehat{B}_s=0\}\cap\{A_s\not=0\})\subset\{\rho_s
\widehat{B}_s=0\} \cap\{A_s\not=0\}.$$
We conclude that $\Phi_s(\mathcal{D})\subset
\{\rho_s\widehat{B}_s=0\}$ holds. Now we are going to show that
the condition $\Phi_s(\mathcal{D}) \subset
\{\rho_s\widehat{B}_s=0\}$ leads to a contradiction. Indeed, we
observe that there exists an irreducible component $\mathcal{D}_0$
of $V_{s+1}$ for which $\mathcal{D}=\A^{(n-s+1)n}\times
\A^{n-s+1}\times \mathcal{D}_0$ holds. Let $x\in\mathcal{D}_0$ be
a nonsingular point of $V_{s+1}$, which is also a nonsingular
point of $V_s$. Hence, for a generic choice of a point
$(\lambda,\gamma)\in\A^{(n-s+1)n}\times \A^{n-s+1}$, the fiber
$W_s:=V_s\cap \{\lambda^*X+\gamma^*=\lambda^*x+\gamma^*\}$ is
unramified (see e.g. \cite[\S 5A]{Mumford95}) and the linear form
$\lambda^{(n-s+1)}X+\gamma_{n-s+1}$ separates the points of $W_s$.
This shows that any point $y\in V_s\cap
\{\lambda^*X+\gamma^*=\lambda^*x+\gamma^*\}$ satisfies the
conditions $D_s(\lambda,\gamma,y)\not=0$ and
$\rho_s(\lambda,\gamma,y)\not=0$. We conclude that the point
$(\lambda,\gamma,\lambda^*x+\gamma^*)$ belongs to the set
$\Phi_s(\mathcal{D})\setminus\{\rho_s\widehat{B}_s=0\}$,
contradicting thus the condition $\Phi_s(\mathcal{D})\subset
\{\rho_s \widehat{B}_s=0\}$. This finishes the proof of our claim.
\spar

From the claim and Lemma \ref{lemma:degree_estimate} we deduce
that the image of the morphism
\[\begin{array}{r}
\!\!\!\!\Psi_s\!:\!(\A^{(n-\!s+\!1)n}\!\!\times\!\A^{n-\!s+\!1}\!\times\!
  V_{s+1})\!\cap\!\{\rho_s\widehat{B}_s\!=0,A_{s+\!1}\!\not=\!0\}\!\!\to
\A^{\!(n-\!s+\!1)n}\!\times\!\A^{n-\!s+\!1}\!\!\times\!\A^{n-\!s-\!1}
\\
(\lambda,\gamma,x)\mapsto\big(\lambda,\gamma,\!Y_1(x),\dots,
Y_{n-\!s-\!1}(x)\big)
\end{array}\]
is contained in a hypersurface of $\A^{(n-s+1)n}\times
\A^{n-s+1}\times\A^{n-s-1}$ of degree at most
$4(n-s+2)^2nd\delta_s^2\delta_{s+1}^2$. Let $\widetilde{B}_s$
denote the defining equation of this hypersurface. \spar

Let
$B_s:=A_sA_{s+1}\rho_s\rho_{s+1}\widehat{B}_s\widehat{B}_{s+1}\widetilde{B}_s$.
Observe that $\deg B_s\le 4(n-s+3)^2nd\delta_s^2\delta_{s+1}^2$
holds. Let $(\lambda,\gamma,P)\in\A^{(n-s+1)n}\times \A^{n-s+1}
\times\A^{n-s}$ be a point satisfying
$B_s(\lambda,\gamma,P)\not=0$. We claim that $(\lambda,\gamma,P)$
satisfies conditions $(i)$, $(ii)$ and $(iii)$ of the statement of
Theorem \ref{theorem:simNoether}. Let $(\lambda^*,\gamma^*)$
denote the first $n-s$ rows of $(\lambda,\gamma)$ and let $P^*$
denote the vector consisting of the first $n-s-1$ coordinates of
$P$. Since $A_s(\lambda,\gamma)A_{s+1} (\lambda^*,\gamma^*)
\not=0$ holds, from Proposition \ref{prop:chow} we conclude that
the mappings $\pi_s:V_s\to\A^{n-s}$ and
$\pi_{s+1}:V_{s+1}\to\A^{n-s-1}$ defined by the linear forms
$\yo{n-s}$ and $\yo{n-s-1}$ are finite morphisms. Since
$A_s(\lambda,\gamma)\not=0$ holds, the condition
$\widehat{B}_s(\lambda,\gamma,P)\not=0$ implies that
$D_s(\lambda,\gamma,x)\not=0$ holds for any $x\in\pi_s^{-1}(P)$.
Therefore, we see that $P$
is a lifting point of the morphism $\pi_s$. A similar argument as
above shows that $P^*$ is a lifting point of the morphism
$\pi_{s+1}$. Finally, the conditions
$\rho_s(\lambda,\gamma,P)\not=0$ and
$\rho_{s+1}(\lambda^*,\gamma^*,P^*)\not=0$ show that $Y_{n-s+1}$
and $Y_{n-s}$ are primitive elements of $\pi_s^{-1}(P)$ and
$\pi_{s+1}^{-1}(P^*)$ respectively. On the other hand, the
conditions $\widetilde{B}_s(\lambda,\gamma,P^*)\not=0$ and
$A_{s+1}(\lambda^*,\gamma^*)\not=0$ imply that
$(\rho_s\widehat{B}_s)\big(\lambda,\gamma,P^*,Y_{n-s}(x)\big)\not=0$
holds for any $x\in\pi_{s+1}^{-1}(P^*)$. Therefore, since
$A_s(\lambda,\gamma)\not=0$ holds, we deduce that
$D_s(\lambda,\gamma,Q)\not=0$ and
$\rho_s(\lambda,\gamma,\pi_s(Q))\not=0$ hold for any point
$Q\in\pi_s^{-1}(P^*,Y_{n-s}(x))$ with $x\in\pi_{s+1}^{-1}(P^*)$.
This shows condition $(iii)$ of the statement of Theorem
\ref{theorem:simNoether}.
\end{proof}

In order to find a rational point of our input variety $V$ we are
going to determine a suitable absolutely irreducible plane
$\fq$--curve of the form $V\cap L$, where $L$ is an
$\fq$--definable affine linear subspace of $\A^n$ of dimension
$r+1$. For this purpose, we are going to find an $\fq$--definable
Noether normalization of $V$, represented by a ($\fq$--definable)
finite linear projection $\pi:V\to\A^{n-r}$, and a lifting point
$P\in\fq^{n-r}$ of $\pi$. Unfortunately, the existence of the
morphism $\pi$ and the point $P$ cannot be guaranteed unless the
number of elements of $\fq$ is high enough. Our next result
exhibits a genericity condition underlying the choice of $\pi$ and
$P$ whose degree depends on $\delta_r:=\deg V_r$, rather than on
$\delta:=\max_{1\le s\le r}\delta_s$.

\begin{corollary}\label{coro:Noether_V_r}
With notations as in Proposition \ref{prop:chow} and Theorem
\ref{theorem:simNoether}, there exists a nonzero polynomial
$\widehat{B}\in\cfq[\Lambda,
\Gamma,\widetilde{Y}_1\klk\widetilde{Y}_{n-r}]$ of degree at most
$(n-r+2)(2nd\delta_r^2-\delta_r)$ such that for any
$(\lambda,\gamma,P)\in\A^{(n-r+1)n}\times \A^{n-r+1}\times
\A^{n-r}$ with $\widehat{B}(\lambda,\gamma,P)\not=0$ the following
conditions are satisfied:

Let $Z:=(Z_1\klk Z_{n-r+1}):=\lambda X+\gamma$. Then the mapping
$\pi:V_r\to\A^{n-r}$ defined by $\pi(x):=\big(Z_1(x)\klk
Z_{n-r}(x)\big)$ is a finite morphism, $P\in\A^{n-r}$ is a lifting
point of $\pi$ and $Z_{n-r+1}$ is a primitive element of
$\pi^{-1}(P)$.
\end{corollary}

\begin{proof}
Let $\widehat{B}:=A_r\rho_r\widehat{B}_r$,  where $A_r$ is the
polynomial of the statement of Proposition \ref{prop:chow} for
$s=r$, $\widehat{B}_r$ is the polynomial of   the proof of Theorem
\ref{theorem:simNoether} with $s=r-1$ and $\rho_r$ is the
discriminant introduced in eq. (\ref{eq:def_disc_chow}) of the
proof of Proposition \ref{prop:chow}. Observe that
$\deg\widehat{B}\le (n-r+2)(2nd\delta_r^2-\delta_r)$ holds. Now,
if $(\lambda,\gamma,P)\in\A^{(n-r+1)n}\times \A^{n-r+1}\times
\A^{n-r}$ is any point for which
$\widehat{B}(\lambda,\gamma,P)\neq 0$ holds, a similar argument as
in the last paragraph of the proof of Theorem
\ref{theorem:simNoether} shows that the linear forms $Z:=\lambda
X+\gamma$ and the point $P$ satisfy the conditions in the
statement of the corollary.
\end{proof}

Combining Theorem \ref{th:Zippel_Schwartz} and Corollary
\ref{coro:Noether_V_r} we conclude that, if $q
>(n-r+2)(2nd\delta_r^2-\delta_r)$ holds, then there  exists an
$\fq$--definable Noether normalization of the variety $V$ and a
lifting point $P\in\fq^{n-r}$ of $\pi$.
%
%
\subsection{A reduction to the bidimensional case}
\label{subsec:bertini} In this section we finish our considerations
about the preparation of the input data
by reducing our problem of computing a rational point
of (the absolutely irreducible $\fq$--variety) $V:=V_r$ to that of
computing a rational point of an absolutely irreducible plane
$\fq$--curve. For this purpose, we have the first Bertini theorem
(see e.g. \cite[\S II.6.1, Theorem 1]{Shafarevich94}), which
asserts that the intersection $V\cap L$ of $V$ with a generic
affine linear subspace $L$ of $\A^n$ of dimension $r+1$ is an
absolutely irreducible plane curve. If $V\cap L$ is an
absolutely irreducible $\fq$--curve, then Weil's estimate (see e.g.
\cite{LiNi83}, \cite{Schmidt76}) assures that we have ``good
probability" of finding a rational point in $V\cap L$. The main
result of this section exhibits an estimate on the degree of the
genericity condition underlying the choice of $L$.

Assume that we are given a point $(\lambda,\gamma,P)\in
\A^{(n-r+1)n} \times\A^{n-r+1}\times \A^{n-r}$ for which
$\widehat{B}(\lambda,\gamma,P)\not=0$ holds, where $\widehat{B}$
is the polynomial of the statement of Corollary \ref{coro:Noether_V_r}.
Let $(Z_1,\dots,Z_{n-r+1})=\lambda X+\gamma$, let $Y_{n-r+2}\klk Y_n$ be
linear forms such that $Z_1,\dots,Z_{n-r+1},Y_{n-r+2}\klk Y_n$ are $\cfq$--linear
independent, and let $P:=(p_1\klk p_{n-r})$. Then the mapping
$\pi:V\to\A^{n-r}$ defined by $\pi(x):=\big(Z_1(x)\klk
Z_{n-r}(x)\big)$ is a finite morphism, and therefore the image
$W:=\widetilde{\pi}(V)$ of $V$ under the mapping
$\widetilde{\pi}:V\to\A^{n-r+1}$ defined by
$\widetilde{\pi}(x):=\big(Z_1(x)\klk Z_{n-r+1}(x)\big)$ is a
hypersurface of $\A^{n-r+1}$. The choice of $Z_1\klk Z_{n-r+1}$
implies that this hypersurface has degree $\delta_r$ and is
defined by a polynomial $q^{(r)}\in\cfq[Z_1\klk Z_{n-r+1}]$ monic
in $Z_{n-r+1}$.

Let $\widetilde{V}:=\{x\in \A^n:(\partial q^{(r)}/\partial
Z_{n-r+1})(Z_1(x),\ldots,Z_{n-r+1}(x))=0\}$ and $\widetilde{W} :=
\{z\in \A^{n-r+1}:(\partial q^{(r)}/\partial Z_{n-r+1})(z)=0\}$.
Our following result shows that the variety $V$ is birationally
equivalent to the hypersurface $W\subset\A^{n-r+1}$.
\begin{lemma}\label{lemma:morph_bir}
$\widetilde{\pi}|_{V\setminus \widetilde{V}}:V\setminus
\widetilde{V}\rightarrow W\setminus \widetilde{W}$ is an
isomorphism of Zariski open sets.
\end{lemma}
\begin{proof}
Let us observe that $\widetilde{\pi}(V\setminus
\widetilde{V})\subset W \setminus \widetilde{W}$. Then
$\widetilde{\pi}|_{V\setminus \widetilde{V}}: V \setminus
\widetilde{V} \to W\setminus\widetilde{W}$ is a well--defined
morphism.

We claim that $\widetilde{\pi}$ is an injective mapping. Indeed,
specializing identity (\ref{equation:chow2}) of the proof of
Proposition \ref{prop:chow} into the values
$\Lambda_{n-r+1,j}:=\lambda_{n-r+1,j}$ ($1\le j\le n$) and
$\Gamma_{n-r+1}=\gamma_{n-r+1}$ we deduce that there exist
polynomials $v_1\klk v_n\in\cfq[Z_1\klk Z_{n-r+1}]$ such that for
$1\le i\le n$ the following identity holds:
\begin{equation}\label{def:v_i}v_i(Z_1\klk Z_{n-r+1})-X_i
\cdot(\partial q^{(r)}/\partial Z_{n-r+1})(Z_1\klk
Z_{n-r+1})\equiv 0 \mbox{ mod $I(V)$ }.
\end{equation}
Let $x:=(x_1\klk x_n),x':=(x_1'\klk x_n')\in V \setminus
\widetilde{V}$ satisfy $\widetilde{\pi}(x)=
\widetilde{\pi}(x')$. We have $Z_i(x)=Z_i(x')$ for $1\le i\le
n-r+1$. Then from identity (\ref{def:v_i}) we conclude that
$x_i=x_i'$ for $1\le i\le n$, which shows our claim.\spar

Now we show that $\widetilde{\pi}|_{V\setminus
\widetilde{V}}:V\setminus \widetilde{V}\to W\setminus
\widetilde{W}$ is a surjective mapping. Let $q_0:=\partial
q^{(r)}/\partial Z_{n-r+1}$. Let be given an arbitrary element
$z:=(z_1\klk z_{n-r+1})$ of $W\setminus \widetilde{W}$, and let
$$x:=\big((v_1/q_0)(z),\ldots,(v_{n}/q_0)(z)\big).$$ We claim that
$x$ belongs to $V\setminus \widetilde{V}$. Indeed, let $F$ be an
arbitrary element of the ideal $I(V)$ and let
$\widetilde{F}:=(q_0(Z_1,\ldots,Z_{n-r+1}))^NF$, where $N:=\deg F$.
Then there exists $G \in\cfq[T_1,\ldots,T_{n+1}]$ such that
$\widetilde{F}=G(q_0X_1,\ldots,q_0X_n,q_0)$ holds. Since $\widetilde{F}
\in I(V)$, for any $z'\in V$ we have $\widetilde{F}(z')=0$, and hence
from identity (\ref{def:v_i}) we conclude that
$G(v_1,\ldots,v_n,q_0)(Z_1(z')\klk Z_{n-r+1}(z'))=0$ holds. This
shows that $q^{(r)}$ divides $\widehat{F}:=G(v_1,\ldots,v_n,q_0)$ in
$\cfq[Z_1\klk Z_{n-r+1}]$ and therefore $\widehat{F}(z)=q_0(z)^NF(x)=0$
holds. Taking into account that $q_0(z)\neq 0$ we conclude that
$F(x)=0$ holds, i.e. $x \in V\setminus \widetilde{V}$.

In order to finish the proof of the surjectivity of
$\widetilde{\pi}$ there remains to prove that
$\widetilde{\pi}(x)=z$ holds. We observe that identity
(\ref{def:v_i}) shows that any $z'\in V$ satisfies
$$Z_{i}(z')q_0\big(Z_1(z')\klk Z_{n-r+1}(z')\big)-\sum_{j=1}^n
\lambda_{i,\,j}\,v_j\big(Z_1(z')\klk Z_{n-r+1}(z')\big)=0$$
for $1\le i\le n-r+1$. Then $q^{(r)}$ divides the polynomial
$Z_iq_0-\sum_{j=1}^n\!\lambda_{i,j}v_j$ in $\cfq[Z_1,\ldots
Z_{n-r+1}]$, which implies $z_i=\sum_{j=1}^n
\lambda_{i,\,j}(v_j/q_0)(z)=\sum_{j=1}^n \lambda_{i,\,j}\,x_j$ for
$1\le i\le n-r+1$. This proves that $\widetilde{\pi}(x)=z$
holds.\spar

Finally we show that $\widetilde{\pi}|_{V\setminus
\widetilde{V}}:V\setminus \widetilde{V}\to W\setminus
\widetilde{W}$ is an isomorphism. Let
\[\begin{array}{ccrcl}
  \phi & : & W\setminus \widetilde{W} & \rightarrow & V\setminus
  \widetilde{V} \\
   &  & z & \mapsto & \big((v_1/q_0)(z),\ldots,
   (v_{n}/q_0)(z)\big).
\end{array}\]
Our previous discussion shows that $\phi$ is a well--defined
morphism. Furthermore, our arguments above show that
$\widetilde{\pi}\circ\phi$ is the identity mapping of
$W\setminus \widetilde{W}$. This finishes the proof of the lemma.
\end{proof}

Let us remark that a similar result for the varieties $V_1\klk
V_{r-1}$ can be easily established following the proof of Lemma
\ref{lemma:morph_bir} {\em mutatis mutandis}. \spar

Now we prove the main result of this section:

\begin{theorem}
\label{theorem:bertini} Let notations and assumptions be as above.
Suppose further that the variety $V:=V_r$ is absolutely
irreducible. Let $\Omega:=(\Omega_1\klk \Omega_{n-r})$ and $T$ be
new indeterminates. Then there exists a nonzero polynomial
$C\in\cfq[\Omega]$ of degree at most $2\delta_r^4$ with the
following property:

\noindent Let $\omega:=(\omega_1\klk\omega_{n-r}) \in\A^{n-r}$
satisfy $C(\omega)\not=0$, and let $L_\omega$ be the
$(r+1)$--dimensional affine linear variety parametrized by
$Z_j=\omega_jT+p_j$ $(1\le j\le n-r)$, $Z_{n-r+1}=Z_{n-r+1}$,
$Y_{n-r+j}=Y_{n-r+j}$
$(2\le j\le r)$. Then $V\cap L_\omega$ is an absolutely
irreducible affine variety of dimension 1.
\end{theorem}
\begin{proof}
Lemma \ref{lemma:morph_bir} shows that $V$ is birational to the
hypersurface $W\subset\A^{n-r+1}$. Since $V$ is absolutely
irreducible, we conclude that $W$ is absolutely irreducible and
therefore $q^{(r)}$ is an absolutely irreducible polynomial.
Following \cite{Kaltofen95a}, let
$\widetilde{q}\in\cfq[\Omega,T][Z_{n-r+1}]$ be the polynomial
$\widetilde{q}:=q^{(r)}\big(\Omega_1T+p_1\klk
\Omega_{n-r}T+p_{n-r},Z_{n-r+1}\big)$.

Since $q^{(r)}$ is a monic element of $\cfq[Z_1\klk
Z_{n-r}][Z_{n-r+1}]$, we easily conclude that $\widetilde{q}$ is a
monic element of $\cfq[\Omega,T][Z_{n-r+1}]$.

We claim that $\widetilde{q}(\Omega,0,Z_{n-r+1})$ is a separable
element of $\cfq[\Omega][Z_{n-r+1}]$. Indeed, we have that
$\widetilde{q}(\Omega,0,Z_{n-r+1})=q^{(r)}(P,Z_{n-r+1})$ holds.
Then the proof of Proposition \ref{prop:chow} shows that the
choice of $P$ implies that the discriminant of the polynomial
$q^{(r)}(P,Z_{n-r+1})$ does not vanish. This means that
$\widetilde{q}(\Omega,0,Z_{n-r+1})$ is a separable element of
$\cfq[\Omega][Z_{n-r+1}]$.

Therefore, applying \cite[Theorem 5]{Kaltofen95a} we conclude that
there exist a polynomial $C\in\cfq[\Omega]$ of degree bounded by
$\frac{3}{2}\delta_r^4-2\delta_r^3+\frac{1}{2}\delta_r^2  \leq
2\delta_r^4$ such that for any $\omega\in\A^{n-r}$ with
$C(\omega)\not=0$, the polynomial
$\widetilde{q}(\omega,T,Z_{n-r+1})$ is absolutely irreducible.
From this we immediately deduce the statement of the theorem.
\end{proof}
%
%
%
\section{The computation of a geometric solution of $V$}
\label{section:computation_geo_sol} Let notations and assumptions
be as in Section \ref{section:preparation}. In this section we
shall exhibit an algorithm which computes a geometric solution of
a ($\K$--definable) lifting fiber $\lfiberr$ of the input variety $V$.

In order to describe this algorithm, we need a simultaneous
Noether normalization of the varieties $V_1\klk V_r$ and lifting
points $P^{(s+1)}\in\A^{n-s-1}$ for $0\le s\le r-1$ such that the
corresponding lifting fiber $V_{P^{(s+1)}}$ has the following
property: for any point $P\in V_{P^{(s+1)}}$, the morphism $\pi_s$
is unramified at $\pi_s(P)$. For this purpose, let
$\Lambda:=(\Lambda_{ij})_{1\le i,j\le n}$ be a matrix of
indeterminates and let $\Gamma:=(\Gamma_1\klk \Gamma_{n})$ be a
vector of indeterminates. Let $X:=(\xon)$ and let
$\widetilde{Y}:=\Lambda X+\Gamma$. Let $B_s\in\cfq[\Lambda,\Gamma,
\widetilde{Y}]$ be the polynomial of the statement of Theorem
\ref{theorem:simNoether} for $1\le s\le r-1$ and let
$B:=\det(\Lambda)\prod_{s=1}^{r-1}B_s$. Observe that $\deg B\le
4n^4d\delta^4$ holds.

Let $\K$ be  a finite field extension of $\fq$ of cardinality
greater than $60n^4d\delta^4$ and let $(\lambda,\gamma,P)$ be a
point randomly chosen in the set ${\K}^{n^2}\times {\K}^n\times
{\K}^{n-r}$. Theorem \ref{th:Zippel_Schwartz} shows that
$B(\lambda,\gamma,P)$ does not vanish with probability at least
$14/15$. From now on, we shall assume that we have chosen
$(\lambda,\gamma,P)\in {\K}^{n^2}\times {\K}^n\times {\K}^{n-r}$
satisfying $B(\lambda,\gamma,P)\not=0$. Let $(Y_1\klk Y_n):=\lambda
X+\gamma$ and $P:=(p_1\klk p_{n-r})$.

From Theorem \ref{theorem:simNoether} we conclude that $Y_1\klk
Y_n$ induce a simultaneous Noether normalization of the varieties
$V_1\klk V_r$ and the point $P^{(s+1)}:=(p_1\klk p_{n-s-1})$
satisfies the condition above for $0\le s\le r-1$. Let us observe
that the fact that the linear forms $\yon$ belong to ${\K}[\xon]$
and $P$ belongs to ${\K}^{n-r}$ immediately implies that the
lifting fiber $\lfibers$ is a $\K$\!--\,variety for $1\le s\le r$.
\spar

The algorithm computing a geometric solution of $\lfiberr$ is a
recursive procedure which proceeds in $r-1$ steps. In the $s$--th
step we compute a geometric solution of the lifting fiber
$V_{P^{(s+1)}}$ from a geometric solution of the lifting fiber
$V_{P^{(s)}}$. Recall that $V_{P^{(s)}}:=\pi_s^{-1}(P^{(s)})=
V_s\cap\{Y_1=p_1\klk Y_{n-s}=p_{n-s}\}$. For this purpose, we
first ``lift" the geometric solution of the fiber $V_{P^{(s)}}$ to
a geometric solution of the affine equidimensional unidimensional
$\K$--variety $\lcurvespone:=V_s\cap\{Y_1=p_1\klk
Y_{n-s-1}=p_{n-s-1}\}$ (see Section \ref{subsec:lifting_to_curve}
below). The variety $\lcurvespone$ is called a {\em lifting curve}.
Then, from this geometric solution we obtain a geometric solution of
the lifting fiber $V_{P^{(s+1)}}=\lcurvespone\cap V(F_{s+1})$. This is
done by computing the minimal equation satisfied by $Y_{n-s+1}$ in
$V_{P^{(s+1)}}$ (see Section \ref{subsec:intersection_step}), from
which we obtain a geometric solution of $V_{P^{(s+1)}}$ by a
suitable effective version of the Shape Lemma (see Section
\ref{subsec:shape_lemma}).
%
%
\subsection{From the lifting fiber \lfibers\ to the lifting curve \lcurvespone}
\label{subsec:lifting_to_curve}
In this section we describe the procedure which computes a
geometric solution of the lifting curve \lcurvespone, from a
geometric solution of the lifting fiber \lfibers.

Let $\pi_s:V_s\to\A^{n-s}$ and $\widetilde{\pi}_{s}:V_s\to
\A^{n-s+1}$ be the linear projection mappings determined by the
linear forms $Y_1\klk
Y_{n-s}$ and $Y_1\klk Y_{n-s+1}$ respectively. From Theorem
\ref{theorem:simNoether} we know that $\pi_s$ is a
finite morphism and $Y_{n-s+1}$ is a primitive element of the
integral ring extension $R_s:=\cfq[\yo{n-s}]\hookrightarrow
\cfq[V_s]$. Furthermore, the minimal polynomial
$q^{(s)}\in\cfq[\yo{n-s+1}]$ of $Y_{n-s+1}$ has degree $\delta_s$
and equals, up to a nonzero element of $\cfq$, the defining
polynomial of the hypersurface $\widetilde{\pi}_{s}(V_s)$. Since
$\widetilde{\pi}_{s}(V_s)$ is a $\K$--hypersurface, we may
assume without loss of generality that $q^{(s)}$ belongs to
${\K}[\yo{n-s+1}]$. This assumption, together with the proof of
Lemma \ref{lemma:morph_bir}, shows that there exists a
geometric solution of $V_s$ consisting of polynomials
$q^{(s)},v_{n-s+2}^{(s)}\klk v_{n}^{(s)}$ of ${\K}[\yo{n-s+1}]$.

We observe that our choice of $P^{(s)}$ implies that the
discriminant of $q^{(s)}$ with respect to $Y_{n-s+1}$ does not
vanish in $P^{(s)}$. Therefore, the above geometric solution of
$V_s$ is compatible with $P^{(s)}$ in the sense of Section
\ref{subsection:lifting_points} and hence the polynomials
$q^{(s)}(P^{(s)},Y_{n-s+1}), v_{n-s+j}^{(s)}(P^{(s)},Y_{n-s+1})$
$(2\le j\le s)$ form a geometric solution of $\lfibers$ with
$Y_{n-s+1}$ as primitive element. We shall assume that we are
given such a geometric solution of $\lfibers$.

Let us observe that \lcurvespone\ can be described as the set of
common zeros of $Y_1-p_1\klk Y_{n-s-1}-p_{n-s-1}, F_1\klk F_s$ or,
equivalently, of $Y_1-p_1\klk
Y_{n-s-1}-p_{n-s-1},F_1(P^{(s+1)},Y_{n-s}\klk Y_n)\klk
F_s(P^{(s+1)},Y_{n-s}\klk Y_n)$. In particular we see
that \lcurvespone\ is a $\K$--variety. In order to find
a geometric solution of \lcurvespone\ we are going to apply the
global Newton--Hensel procedure of \cite{GiLeSa01}. For this
purpose, we need the following result.
\begin{lemma}\label{lemma:sucesionregular}
$F_1(P^{(s+\!1)}\!,Y_{n-\!s},\dots,\! Y_n),\dots,\!
F_s(P^{(s+\!1)}\!,Y_{n-\!s},\dots,\! Y_n)$ generate a radical
ideal of ${\K}[Y_{n-s}\klk Y_n]$ and form a regular sequence of
${\K}[Y_{n-s}\klk Y_n]$, and
$\lcurvespone$ has degree $\delta_s$.

\end{lemma}
\begin{proof}
We first show that  $F_1(P^{(s+\!1)},Y_{n-s}\ldots Y_n ),\ldots,
F_s(P^{(s+\!1)}, Y_{n-s},\ldots, Y_n)$ form a regular sequence.
Let $L_{s+1}\subset\A^n$ be the affine linear variety
$L_{s+1}:=\{Y_1=p_1\klk Y_{n-s+1}=p_{n-s+1}\}$. Observe that, for
$1\le i\le s$, the discriminant $\rho^{(i)}$ of the polynomial
$q^{(i)}$ with respect to $Y_{n-i+\!1}$ does not vanish in
$P^{(i)}:=(p_1\klk p_{n-i})$. This implies that
$q^{(i)}(P^{(s+1)},Y_{n-s}\klk\! Y_{n-i+1})$ is a separable
polynomial of \linebreak ${\K}[Y_{n-s}\klk Y_{n-i+1}]$ for $1\le
i\le s$. Hence, $V_i\cap L_{s+1}\cap \{\partial q^{(i)}/\partial
Y_{n-i+1}\not=0\}$ is a nonempty Zariski--dense open set of
$V_i\cap L_{s+1}$. Following Lemma \ref{lemma:morph_bir}, let
$\widetilde{V}_i:=\{x\in\A^n:(\partial q^{(i)}/\partial
Y_{n-i+1})(x)\not=0\}$, $W_i:=\widetilde{\pi}(V_i)$ and
$\widetilde{W}_i:=\widetilde{\pi}_i(\widetilde{V}_i)$ for $1\le i
\le s$. We have $\widetilde{\pi}_i((V_i\setminus
\widetilde{V}_i)\cap L_{s+1})=(W_i\setminus\widetilde{W}_i)\cap
L_{s+1}$. Therefore, since $(W_i\setminus\widetilde{W}_i)\cap
L_{s+1}$ has dimension $s+1-i$, and taking into account that
$\widetilde{\pi}_i|_{V_i\setminus
\widetilde{V}_i}:V_i\setminus\widetilde{V}_i\to
W_i\setminus\widetilde{W}_i$ is an isomorphism of locally closed
sets, we conclude that $(V_i\setminus\widetilde{V}_i)\cap L_{s+1}$
has dimension $s+1-i$, which implies that $V_i\cap L_{s+1}$ has
dimension $s+1-i$ for $1\le i\le s$. This proves that
$F_1(P^{(s+1)},Y_{n-s}\klk Y_n )\klk F_s(P^{(s+1)}, Y_{n-s}\klk
Y_n)$ form a regular sequence of $\K[Y_{n-s}\klk Y_n]$.\spar

Now we prove that $\deg\lcurvespone=\delta_s$ holds. Observe that
our previous argumentation shows that $\lcurvespone=V_s\cap
L_{s+1}$ is an equidimensional variety of dimension 1 which, by
the B\'ezout inequality (\ref{equation:Bezout}), satisfies the
degree estimate $\deg\lcurvespone\le\delta_s$. On the other hand,
since $\pi_s$ is a finite morphism we conclude that the
restriction mapping $\pi_s|_{\lcurvespone}:\lcurvespone\to
L_{s+1}$ is also a finite morphism. Furthermore, our choice of
$P^{(s)}$ implies that
$\#(\pi_s|_{\lcurvespone})^{-1}(P^{(s)})=\#\,\pi_s^{-1} (P^{(s)})=
\delta_s$ holds. Then
$$\delta_s=\#\pi_s^{-1}(P^{(s)})=\#(\lcurvespone\cap
\{Y_{n-s}=p_{n-s}\})\le\deg\lcurvespone\le\delta_s,$$ which proves
our second assertion. \spar

There remains to prove that $F_1(P^{(s+\!1)}\!,Y_{n-\!s},\dots,\!
Y_n),\dots,\! F_s(P^{(s+\!1)}\!,Y_{n-\!s},\dots,\! Y_n)$ generate
a radical ideal of ${\K}[Y_{n-s}\klk Y_n]$. Let
$\pi_{s+1}^*:\lcurvespone\to\A^{s+1}$ be the projection defined by
$Y_{n-s}\klk Y_n$. Observe that $\lcurvesponestar:=
\pi_{s+1}^*(\lcurvespone)$ is the subvariety of $\A^{s+1}$ defined by
$F_1(P^{(s+1)},Y_{n-s},\dots, Y_n),\dots,
F_s(P^{(s+1)},Y_{n-s},\dots, Y_n)$ and is isomorphic to
$\lcurvespone$. Since
$\lcurvesponestar\cap\{Y_{n-s}=p_{n-s}\}=\pi_{s+1}^*(\lfibers)$
holds, from Lemma \ref{lemma:equiv_ramif} we conclude that the
Jacobian determinant
$$J_F(P^{(s+1)},Y_{n-s}\klk Y_n):=\det\big(\partial
F_i(P^{(s+1)},Y_{n-s+1}\klk Y_n)/\partial Y_{n-s+j}\big)_{1\le
i,j\le s}$$
does not vanish in any point
$P\in\lcurvesponestar\cap\{Y_{n-s}=p_{n-s}\}$.
Furthermore, the identity $\#(\lcurvesponestar \cap
\{Y_{n-s}=P_{n-s}\})=\delta_s=\deg \lcurvesponestar$ shows that
the affine linear variety $\{Y_{n-s}=p_{n-s}\}$ meets every
irreducible component of $\lcurvesponestar$.
This proves that the coordinate
function of \lcurvespone\ defined by
$J_F(P^{(s+1)},Y_{n-s}\klk Y_n)$
is not a zero divisor of $\cfq[\lcurvespone]$. Hence,
from \cite[Theorem 18.15]{Eisenbud95} we conclude that the ideal
generated by $F_1(P^{(s+1)}, Y_{n-s}\klk Y_n)\klk F_s(P^{(s+1)},
Y_{n-s}\klk Y_n)$ is radical.
\end{proof}
Using the notations of the proof of Lemma
\ref{lemma:sucesionregular}, let $\pi_{s+1}^*:\lcurvespone\to
\A^{s+1}$ be the projection mapping defined by $Y_{n-s}\klk Y_n$,
and let $\lcurvesponestar:= \pi_{s+1}^*(\lcurvespone)$ and
$\lfibersstar:= \pi_{s+1}^*(\lfibers)$. Then the projection
$\widehat{\pi}_{s+1}:\lcurvesponestar\to\A^1$ induced by $Y_{n-s}$ is
a finite morphism of degree $\delta_s$, whose fiber
$\widehat{\pi}_{s+1}^{-1}(p_{n-s})=\lfibersstar$ is unramified.
Furthermore, the polynomials
$q^{(s)}(P^{(s)},Y_{n-s+1}),v_{n-s+j}^{(s)}(P^{(s)},Y_{n-s+1})$
$(2\le j\le s)$, introduced before the statement of Lemma
\ref{lemma:sucesionregular}, form a geometric solution of
$\lfibersstar$. Under these conditions, applying the Global Newton
algorithm of \cite[II.4]{GiLeSa01} we conclude that there exists a
computation tree $\beta$ in $\K$ which computes a geometric
solution of $\lcurvesponestar$, which is also a geometric solution
of $\lcurvespone$. Let us observe that the fact that the input
geometric solution of $\lfibersstar$ consists of univariate
polynomials with coefficients in $\K$ implies that the output
geometric solution of $\lcurvespone$ also consists of polynomials
with coefficients in $\K$.

The evaluation of the computation tree $\beta$ requires
$O\big((n\mathcal{T}+n^5) \mathcal{U}(\delta_s)^2\big)$ arithmetic
operations in ${\K}$, using at most
$O\big((\mathcal{S}+n)\delta_s^2\big)$ arithmetic registers. Here
$\mathcal{S}$ and $\mathcal{T}$ denote the space and time of the
input \slp\ representing the polynomials $F_1\klk F_s$ and
$\mathcal{U}(m)$ denotes the quantity
$\mathcal{U}(m):=m\log^2m\log\log m$ for any $m\in\N$. Let us
remark that the asymptotic estimate $O\big(\mathcal{U}(m)\big)$
represents the bit--complexity of certain basic operations (such as
addition, multiplication, division and gcd) with integers of
bit--size $m$, and the number of arithmetic operations in a given
domain $R$ necessary to compute the multiplication, division,
resultant, gcd and interpolation of univariate polynomials of
$R[T]$ of degree at most $m$ (cf. \cite{GaGe99}, \cite{BiPa94}).
In particular, an arithmetic operation in a finite field  ${\K}$
of cardinality $\# {\K} $ can be (deterministically) performed with
$O(\mathcal{U}(\log \# {\K}))$ bit operations, using space $
O(\log \# {\K}).$ Our assumptions on ${\K}$
imply $\log \# {\K} \le O(\log (q \delta))$.\spar

From the above considerations we easily deduce the following
result:
\begin{proposition}\label{prop:lifting_to_curve}
There exists a (deterministic) Turing machine $M$ which has as
input
\begin{itemize}
\item
a \slp\ using space $\mathcal{S}$ and time $\mathcal{T}$ which
represents the polynomials $F_1\klk F_s$,
\item
the dense representation of elements of ${\K}[Y_{n-s+1}]$ which
form a geometric solution of $\lfibers$,
\end{itemize}
and outputs the dense representation of polynomials of
${\K}[Y_{n-s},Y_{n-s+1}]$ which form a geometric solution of
\lcurvespone. The Turing machine $M$ runs in space
$O\big((\mathcal{S}+n)\delta_s^2\log (q\delta)\big)$ and time
$O\big((n\mathcal{T}+n^5) \mathcal{U}(\delta_s)^2 \mathcal{U}(\log
(q\delta))\big)$.
\end{proposition}
%
%
\subsection{Computing a hypersurface birational to $V_{P^{(s+1)}}$ }
\label{subsec:intersection_step}
The purpose of this section is to exhibit an algorithm which
computes the minimal equation satisfied by the coordinate function
induced by a linear form $\Llambda:=Y_{n-s}+\lambda Y_{n-s+1}$
in $\cfq[\lfiberspone]$, for a suitable choice of $\lambda\in K$.

With the notations of the previous section, let
$\pi_{s+1}^*:\lcurvespone\to\A^{s+1}$ be the projection defined by
$Y_{n-s}\klk Y_n$, let $\lcurvesponestar:=
\pi_{s+1}^*(\lcurvespone)$ and let $\lfibersponestar:=
\pi_{s+1}^*(\lfiberspone)$. For a given $\lambda \in K$, let
$\mathcal{L}_\lambda\in {\K}[Y_{n-s},Y_{n-s+1}]$ denote the linear
form $\mathcal{L}_{\lambda}:= Y_{n-s}+\lambda Y_{n-s+1}$, and let
$\widehat{\pi}_{s+1,\lambda}:\lcurvesponestar\to\A^{1}$ be the
projection morphism defined by $\widehat{\pi}_{s+1,\lambda}(x)
:=\mathcal{L}_\lambda (x)$. Our next result yields a sufficient
(and consistent) condition on $\lambda$, which assures that
replacing the variable $Y_{n-s}$ by $\mathcal{L}_{\lambda}$ does
not change the situation obtained after the preprocessing of
Section \ref{subsec:main_preparation}, namely
$\widehat{\pi}_{s+1,\lambda}$ is a finite morphism and any element
of the set $\widehat{\pi}_{s+1,\lambda}(\lfibersponestar)$ is an
unramified point of $\widehat{\pi}_{s+1,\lambda}$.
\begin{lemma}\label{lemma:moving_projection}
Let $\Lambda$ be an indeterminate. There exists a nonzero
polynomial $E_s\in\cfq[\Lambda]$ of degree at most $4\delta^3$,
with the following property:\sspar

\noindent Let $\lambda\in\A^1$ be any point with
$E_s(\lambda)\not=0$, and let $\mathcal{L}_\lambda:=Y_{n-s}+\lambda
Y_{n-s+1}$. Then the following conditions are satisfied:
\begin{itemize}
\item[$(i)$] The projection mapping
$\widehat{\pi}_{s+1,\lambda}:\lcurvespone^*\to\A^1$ determined by
$\mathcal{L}_\lambda$ is a finite morphism.
\item[$(ii)$] $\mathcal{L}_\lambda$ separates the points
of the lifting fiber $\lfibersponestar$.
\item[$(iii)$] Every element of $\widehat{\pi}_{s+1,\lambda}(\lfiberspone^*)$
is a lifting point of $\widehat{\pi}_{s+1,\lambda}$.
\end{itemize}
\end{lemma}
\begin{proof}
By the choice of the linear forms $Y_1\klk Y_n$ and the point
$P^{(s+1)}$ we have that $Y_{n-s+1}$ induces a primitive element
of the integral ring extension $\cfq[Y_{n-s}]\hookrightarrow
\cfq[\lcurvespone]$, whose minimal polynomial is
$q^{(s)}(P^{(s+1)},Y_{n-s},Y_{n-s+1})$. Furthermore,
$\cfq[\lcurvespone]$ is a free $\cfq[Y_{n-s}]$--module of rank
$\delta_s$.

First we determine a genericity condition which assures that
condition $(i)$ of the statement of Lemma
\ref{lemma:moving_projection} is satisfied. Let $\Lambda$ be a new
indeterminate, let $\LLambda:= Y_{n-s}+\Lambda Y_{n-s+1}$, and let
$q_{\Lambda}^{(s)}$ be the following element of
${\K}[\Lambda,Y_1\klk Y_{n-s-1},\LLambda,Y_{n-s+1}]$:
$$q_{\Lambda}^{(s)}:= q^{(s)}(Y_1\klk Y_{n-s-1},\LLambda-\Lambda
Y_{n-s+1},Y_{n-s+1}).$$
Since $q^{(s)}$ has (total) degree $\delta_s$ and the expression
$\LLambda-\Lambda Y_{n-s+1}$ is linear in the variables $\LLambda$
and $Y_{n-s+1}$, we conclude that
$\deg_{\LLambda,Y_{n-s+1}}q_{\Lambda}^{(s)}\le\delta_s$ holds. On
the other hand it is clear that $q_{\Lambda}^{(s)}$ has degree at
most $\delta_s$ in $\Lambda$. Therefore, we may express
$q_{\Lambda}^{(s)}(P^{(s+1)},\Lambda,\LLambda,Y_{n-s+1})$ in the
following way:
$$q_{\Lambda}^{(s)}(P^{(s+1)},\Lambda,\LLambda,Y_{n-s+1})=
a_{\delta_s}(\Lambda)Y_{n-s+1}^{\delta_s}+
a_{\delta_s-1}(\Lambda,\LLambda)Y_{n-s+1}^{\delta_s}\plp
a_0(\Lambda,\LLambda),$$
where $a_{\delta_s},\dots\!, a_0\!\in\!{\K}[\Lambda,\!\LLambda]$
have degree at most $\delta_s$. Since
$q_{\Lambda}^{(s)}\!(P^{(s+\!1)}\!,0,Y_{\!n-s},Y_{\!n-s+\!1})$ $=
q^{(s)}(P^{(s+1)},Y_{n-s},Y_{n-s+1})$ holds and
$q^{(s)}(P^{(s+1)},Y_{n-s},Y_{n-s+1})$ is a monic element of
${\K}[Y_{n-s}][Y_{n-s+1}]$ of degree $\delta_s$ in $Y_{n-s+1}$, we
conclude that the leading coefficient $a_{\delta_s}$ is a nonzero
element of ${\K}[\Lambda]$ (of degree at most $\delta_s$). We
shall prove below that for any $\lambda$ with
$a_{\delta_s}(\lambda)\not=0$ condition $(i)$ holds. \spar

Now we consider condition $(ii)$ of the statement of Lemma
\ref{lemma:moving_projection}. Let $\lfibersponestar:=\{Q_1\klk
Q_{\delta_{s+1}}\}$, and let us consider the following polynomial:
$$E_{s,1}(\Lambda)=\!\!\!\!\!\!\!\!\prod_{\quad 1\le j<k\le
\delta_{s+1}}\!\!\!\!\! \!\!\!\!\big( \LLambda(Q_j)- \LLambda(
Q_k)\big).$$
Observe that $\LLambda\!(Q_j)-\!\LLambda\!(
Q_k)\!=\!Y_{n-\!s}(Q_j)-\!Y_{n-\!s}(Q_k)+
\Lambda\big(Y_{n-\!s+\!1}(Q_j)-\!Y_{n-\!s+\!1}(Q_k)\big)$ holds
for $1\le j<k\le \delta_{s+1}$. Therefore, since $Y_{n-s}$
separates the points of the lifting fiber \lfibersponestar, we
conclude that $E_{s,1}$ is a nonzero element of $\cfq[\Lambda]$ of
degree at most $\delta_{s+1}^2$. We shall show below that for any
$\lambda$ with $E_{s,1}(\lambda)\not=0$ condition $(ii)$ holds. \spar

Finally, we analyze condition $(iii)$ of the statement of Lemma
\ref{lemma:moving_projection}. Let
$\widehat{\pi}_{s+1,\Lambda}:\A^1\times\lfiberspone\to\A^2$ be the mapping
defined by $\widehat{\pi}_{s+1,\Lambda}(\lambda,x):=\big(\lambda,
\mathcal{L}_\lambda(x)\big)$. Observe that the image of
$\widehat{\pi}_{s+1,\Lambda}$ is a $\K$--hypersurface of $\A^2$ of degree
$\delta_{s+1}$, which is defined by the polynomial
$q_{\LLambda}^{(s+1)}(P^{(s+1)},\Lambda,\LLambda):=\prod_{1\le
j\le \delta_{s+1}}(\LLambda-\LLambda(Q_j)) \in
{\K}[\Lambda,\LLambda]$. We claim that
$q_{\LLambda}^{(s+1)}(P^{(s+1)},\Lambda,\LLambda)$ and the
discriminant $\rho_{\Lambda}^{(s)}(P^{(s+1)},\Lambda,\LLambda)\in
{\K}[\Lambda, \LLambda]$ of the polynomial
$q_{\Lambda}^{(s)}(P^{(s+1)}, \Lambda,\LLambda, Y_{n-s+1})$
introduced above have no nontrivial common factors in
${\K}(\Lambda)[\LLambda]$. Arguing by contradiction, suppose that
there exists a nontrivial common factor $\widetilde{h}\in
{\K}(\Lambda)[\LLambda]$ of
$q_{\LLambda}^{(s+1)}(P^{(s+1)},\Lambda,\LLambda)$ and
$\rho_{\Lambda}^{(s)}(P^{(s+1)},\Lambda,\LLambda)$. Since
$q_{\LLambda}^{(s+1)}(P^{(s+1)},\Lambda,\LLambda)$ is a monic
element of ${\K}[\Lambda][\LLambda]$, we deduce that there exists a
common factor $h\in {\K}[\Lambda,\LLambda] \setminus
{\K}[\Lambda]$ not divisible by $\Lambda$. Taking into account
that $q_{\LLambda}^{(s+1)}
(P^{(s+1)},0,Y_{n-s})=q^{(s+1)}(P^{(s+1)},Y_{n-s})$ holds and
$\rho_{\Lambda}^{(s)}(P^{(s+\!1)}\!,0,\!Y_{n-s})$ is the
discriminant $\rho^{(s)}(P^{(s+\!1)}\!,Y_{\!n-s})$ of $q^{(s)}
(P^{(s+\!1)}\!,Y_{\!n-s}, Y_{\!n-s+1})$ with respect to
$Y_{n-s+1}$, we conclude that $h(0,Y_{n-s})$ is a nontrivial
common factor of $\rho^{(s)}(P^{(s+1)},Y_{n-s})$ and
$q^{(s+1)}(P^{(s+1)},Y_{n-s})$. Let $\alpha\in\cfq$ be a root of
$h(0, Y_{n-s})$ and let $Q$ be a point of $\lfiberspone$ for which
$\alpha=Y_{n-s}(Q)$ holds. Then $(p_1\klk
p_{n-s-1},\alpha)=\pi_s(Q)$, and $q^{(s)}(\pi_s(Q),Y_{n-s+1})$ has
less than $\delta_s$ roots. We conclude that either
$\pi_s(Q)$ is not a lifting point of $\pi_s$ or
$Y_{n-s+1}$ is not a primitive element of $\pi_s^{-1}(\pi_s(Q))$,
contradicting thus condition $(iii)$ of Theorem
\ref{theorem:simNoether}. This proves our claim.

From our claim we see that the resultant $E_{s,2}\in{\K}[\Lambda]$
of $q_{\LLambda}^{(s+1)}(P^{(s+1)},\Lambda,\LLambda)$ and
$\rho_{\Lambda}^{(s)}(P^{(s+1)},\Lambda,\LLambda)$ with respect to
the variable $\LLambda$ is a nonzero element of $\cfq[\Lambda]$ of
degree at most $2(2\delta_s-1)\delta_s\delta_{s+1}$. The
nonvanishing of $E_{s,2}$ is the genericity condition which
implies condition $(iii)$, as will be shown below. \spar

Let $E_s:=a_{\delta_s}E_{s,1}E_{s,2}\in\cfq[\Lambda]$. Observe
that $\deg E_s\le 4\delta^3$ holds. Let $\lambda\in\A^1$ satisfy
$E_s(\lambda)\not=0$ and let $\Llambda:=Y_{n-s}+\lambda
Y_{n-s+1}$. We claim that conditions $(i)$, $(ii)$ and $(iii)$ of
the statement of Lemma \ref{lemma:moving_projection} hold.

Let $\ell_\lambda$, $y_{n-s}$ and $y_{n-s+1}$ denote the
coordinate functions of $\cfq[\lcurvespone]$ induced by
$\mathcal{L}_\lambda:=Y_{n-s}+\lambda Y_{n-s+1}$, $Y_{n-s}$ and
$Y_{n-s+1}$ respectively. We have $\ell_\lambda=y_{n-s}+ \lambda
y_{n-s+1}$. From the identity $q^{(s)}(P^{(s+1)},y_{n-s},
y_{n-s+1})=0$ we deduce that
$q_{\Lambda}^{(s)}(P^{(s+\!1)},\lambda,\ell_\lambda,
y_{n-s+\!1})=0$ holds. Let
$q_{\lambda}^{(s)}\!:=q_{\Lambda}^{(s)}(\lambda, Y_1\klk
Y_{n-s-1},\Llambda,Y_{n-s+1})$. Since
$a_{\delta_s}(\lambda)\not=0$ holds, we see that
$q_{\lambda}^{(s)}(P^{(s+1)},\Llambda,Y_{n-s+1})$ is a monic (up
to a nonzero element of $\cfq$) element of
$\cfq[\Llambda][Y_{n-s+1}]$, which represents an integral
dependence equation over $\cfq[\Llambda]$ for the coordinate
function $y_{n-s+1}$. Assuming without loss of generality that
$\lambda\not=0$ holds, we see that
$\widehat{\pi}_{s+1,\lambda}:\lcurvesponestar\to\A^1$ is a dominant
mapping, for
otherwise $\widehat{\pi}_{s+1}:\lcurvesponestar\to\A^1$ would not be dominant.
Therefore, we conclude that $\cfq[\Llambda]
\hookrightarrow\cfq[\llambda,y_{n-s+1}]$ is an integral ring
extension, which, combined with the fact that
$\cfq[\llambda,y_{n-s+1}] \to \cfq[\lcurvesponestar]$ is an integral
ring extension, implies that
$\cfq[\Llambda]\hookrightarrow\cfq[\lcurvesponestar]$ is an integral
extension. This proves that $\widehat{\pi}_{s+1,\lambda}$ is a finite
morphism and shows that condition $(i)$ holds.

Next, taking into account that $E_{s,1}(\lambda)=\prod_{1\le
i<j\le \delta_{s+1}}\big(\Llambda(Q_i)-\Llambda(Q_j)\big)\not=0$ holds,
we conclude that $\Llambda$ separates the points of the fiber
$\lfibersponestar$. This shows that condition $(ii)$ holds.

Finally, let $Q$ be an arbitrary point of $\lfibersponestar$. Since
$E_{s,2}(\lambda)\not=0$ holds, we have that the discriminant
$\rho_{\lambda}^{(s)}(P^{(s+1)},\Llambda)$ of the polynomial
$q_{\lambda}^{(s)}(P^{(s+1)},\Llambda,Y_{n-s+1})$ with respect to
$Y_{n-s+1}$ does not vanish in $\Llambda(Q)$. This implies that the
polynomial $q_{\lambda}^{(s)}
(P^{(s+1)},\Llambda(Q),\!Y_{n-\!s+\!1})$ has $\delta_s$ distinct
roots in $\cfq$, which proves that the fiber
$\widehat{\pi}_{s+1,\lambda}^{-1} (\Llambda(Q))$ has $\delta_s$ distinct
points, i.e. is unramified. This shows that condition $(iii)$ holds
and finishes the proof of the lemma.\end{proof}

Since the cardinality of the field ${\K}$ is greater than
$60n^4d\delta^4$, from Theorem \ref{th:Zippel_Schwartz} we see
that, for a randomly chosen value $\lambda \in {\K}$, the condition
$E_s(\lambda)\not=0$ holds with probability at least $1-1/60n^4$.
Assume that we are given such a value $\lambda\in {\K}$ and let
$\Llambda:=Y_{n-s}+\lambda Y_{n-s+1}$. We are going to exhibit an
algorithm that computes the minimal equation of the coordinate
function of $\lfiberspone$ induced by {\Llambda}.

Let $(\partial q_\lambda^{(s)}\!/\partial
Y_{\!n-\!s+\!1})^{-1}(P^{(s+1)}, \Llambda,Y_{\!n-\!s+\!1})$ be the
monic element of ${\K}(\Llambda)[Y_{\!n-\!s+\!1}]$ of degree at
most $\delta_s-1$ that is the inverse of $(\partial
q_\lambda^{(s)}\!/\partial
Y_{\!n-\!s+\!1})(P^{(s+1)}\!,\Llambda,Y_{\!n-\!s+\!1})$ modulo
$q_\lambda^{(s)}(P^{(s+1)}, \Llambda,Y_{n-s+1})$, and let
$w_{\!n-\!s+\!j}^{(s)}(P^{(s+1)},\Llambda,Y_{\!n-\!s+\!1})\in {\K}
(\Llambda) [Y_{\!n-\!s+\!1}]$ be the remainder of the product
$v_{\!n-\!s+\!j}^{(s)}(P^{(s+1)}\!,\Llambda-\lambda
Y_{n-s+1},Y_{n-s+1})(\partial q_\lambda^{(s)}\!/\partial
Y_{\!n-\!s+\!1})^{-1}$ $(P^{(s+1)},\Llambda,Y_{\!n-\!s+\!1})$
modulo $q_\lambda^{(s)}(P^{(s+1)},\Llambda,Y_{\!n-\!s+\!1})$ for
$2\le j\le s$. Finally, let

$$f_{\!s+\!1}\!:=\!F_{\!s+\!1}\!\Big(\!P^{(\!s+\!1)}\!,\Llambda,\!
Y_{\!n\!-\!s+\!1},w_{n\!-\!s+2}^{(\!s)}(P^{(\!s+\!1)}\!,
\Llambda,\!Y_{\!n\!-\!s+\!1}),\dots,
w_n^{(\!s)}(P^{(\!s+\!1)}\!,\Llambda,\!Y_{\!n\!-\!s+\!1})\!\Big),$$
\begin{equation}\label{eq:expression_g}
g_{s+1}:=\!\mbox{Res}_{Y_{n-s+1}}\Big(q^{(s)}(P^{(s+1)}
,\Llambda,Y_{n-s+1}),f_{s+1}\Big),\end{equation}
where $\mathrm{Res}_{Y_{n-s+1}}(f,g)$ denotes the resultant of
$f$ and $g$ with respect to $Y_{n-s+1}$.

Let us observe that $f_{s+1}\in {\K}(\Llambda)[Y_{n-s+1}]$ has
degree at most $d\delta_s$ (in $Y_{n-s+1}$), and the denominators
of its coefficients are divisors of a polynomial of
${\K}[\Llambda]$ of degree bounded by $(2\delta_s-1)\delta_s$. On the
other hand, from \cite[Corollary 2]{GiLeSa01} it follows that
$g_{s+1}$ is an element of ${\K}[\Llambda]$ of degree bounded by
$d\delta_s$. Our next result shows that the minimal equation of
$\Llambda$ in ${\K}[\lfiberspone]$ can be efficiently computed.
\begin{proposition} \label{prop:minimal}
 There exists a probabilistic Turing machine $M$ which has as
input
\begin{itemize}
\item a \slp\ using space $\mathcal{S}$ and time $\mathcal{T}$
which represents the polynomial $F_{s+1}$,
\item the dense representation of elements of
${\K}[Y_{n-s},Y_{n-s+1}]$ which form a geometric solution of
$\lcurvespone$, as computed in Proposition
\ref{prop:lifting_to_curve},
\item a value $\lambda\in {\K}$ satisfying the conditions of Lemma
\ref{lemma:moving_projection},
\end{itemize}
and outputs the dense representation of the minimal polynomial
$q_{\Llambda}^{(s+1)}(P^{(s+1)},\Llambda)$ $\in {\K}[\Llambda]$ of
the coordinate function of $\lfiberspone$ induced by \Llambda. The
Turing machine $M$ runs in space
$O\big((\mathcal{S}+d)\delta_s^2\log (q\delta)\big)$ and time
$O\big((\mathcal{T}+n)\mathcal{U}(d\delta_s)
\mathcal{U}(\delta_s)\mathcal{U}(\log (q\delta))\big)$ and outputs
the right result with probability at least $1-1/45n^3$.
\end{proposition}
\begin{proof}
Let $\lambda\in {\K}$ satisfy the conditions of Lemma
\ref{lemma:moving_projection}. Then \cite[Lemma 8]{HeMaWa01} shows
that the following identity holds:
$$q^{(s+1)}_{\Llambda}(P^{(s+1)},\Llambda)=\frac{g_{s+1}}{\gcd
(g_{s+1},g_{s+1}'\big)}\,.$$
Therefore, the computation of $q^{(s+1)}_{\Llambda}(P^{(s+1)},
\Llambda)$ can be efficiently reduced to that of the polynomial
$g_{s+1}$ of (\ref{eq:expression_g}). The latter may be defined as
the resultant with respect to the variable $Y_{n-s+1}$ of two
elements of ${\K}(Y_{n-s})[Y_{n-s+1}]$ of degrees bounded by
$\delta_s$ and $\delta_s-1$, namely
$q^{(s)}(P^{(s+1)},Y_{n-s},Y_{n-s+1})$ and the remainder of
$f_{s+1}$ modulo $q^{(s)}(P^{(s+1)},Y_{n-s},Y_{n-s+1})$. Following
\cite[Corollary 11.16]{GaGe99}, such resultant can be computed
using the Extended Euclidean Algorithm (EEA for short) in
${\K}(\Llambda)[Y_{n-s+1}]$, which requires
$O\big(\mathcal{U}(\delta_s)\big)$ arithmetic operations in
${\K}(\Llambda)$ storing at most $O(\delta_s)$ elements of
${\K}(\Llambda)$. Furthermore, the computation of $f_{s+1}$
requires the (modular) inversion of $(\partial
q_\lambda^{(s)}/\partial Y_{n-s+1})^{-1}(P^{(s+1)},
\Llambda,Y_{n-s+1})$, which can also be computed applying the EEA
in ${\K}(\Llambda)[Y_{n-s+1}]$ to the polynomials
$q_\lambda^{(s)}(P^{(s+1)},\Llambda,Y_{n-s+1})$ and $(\partial
q_\lambda^{(s)}/\partial
Y_{n-s+1})(P^{(s+1)},\Llambda,Y_{n-s+1})$.

In order to compute the dense representation of the polynomial
$g_{s+1}$, we shall perform the EEA over a ring of power series
${\K}[\![\Llambda-\alpha]\!]$ for some ``lucky" point $\alpha \in
{\K}$. Therefore, we have to determine a value $\alpha\in {\K}$
such that all the elements of ${\K}[\Llambda]$ which are inverted
during the execution of the EEA are invertible elements of the
ring ${\K}[\![\Llambda-\alpha]\!]$. Further, in order to make our
algorithm ``effective", during its execution we shall compute
suitable approximations in ${\K}[\Llambda]$ of the intermediate
results of our computations, which are obtained by truncating the
power series of ${\K}[\![\Llambda-\alpha]\!]$ that constitute
these intermediate results. Therefore, we have to determine the
degree of precision of the truncated power series required to
output the right results.

In order to determine the value $\alpha\in {\K}$, we observe that,
similarly to the proof of \cite[Theorem 6.52]{GaGe99}, one deduces
that all the denominators of the elements of ${\K}(\Llambda)$
arising during the application of the EEA to
$q_\lambda^{(s)}(P^{(s+1)},\Llambda,Y_{n-s+1})$ and $f_{s+1}$ are
divisors of at most $\delta_{s}+1$ polynomials of ${\K}[Y_{n-s}]$
of degree bounded by $(d\delta_s+\delta_s)(2\delta_s-1)\delta_s$.
On the other hand, the denominators arising during the application
of the EEA to $q_\lambda^{(s)}(P^{(s+1)},\Llambda,Y_{n-s+1})$ and
$(\partial q_\lambda^{(s)}/\partial
Y_{n-s+1})(P^{(s+1)},\Llambda,Y_{n-s+1})$ are divisors of at most
$\delta_{s}+1$ polynomials of ${\K}[Y_{n-s}]$ of degree at most
$(2\delta_s-1)\delta_s$. Hence the product of all the denominators
arising during the two applications of the EEA has degree at most
$(d\delta_s+\delta_s+1)(2\delta_s-1)\delta_s(\delta_s+1)\le
4d\delta_s^4$. Since $\# {\K} >60n^4d\delta^4$ holds, from Theorem
\ref{th:Zippel_Schwartz} we conclude that there exists $\alpha\in
K $ that does not vanish any denominator arising as an
intermediate results of the EEA. Furthermore, the probability of
finding such $\alpha$ by a random choice in $\K$ is at least
$1-1/45n^3$.

On the other hand, since the output of our algorithm is a
polynomial of degree at most $d\delta_s$, computing all the power
series which arise as intermediate results up to order
$d\delta_s+1$ allows us to output the right result.

Our algorithm computing $g_{s+1}$ inverts $(\partial
q_\lambda^{(s)}/\partial
Y_{n-s+1})(P^{(s+1)}\!,\Llambda,Y_{n-s+1})$ modulo
$q_\lambda^{(s)}(P^{(s+1)},\Llambda,Y_{n-s+1})$, computes
$w_{n-s+j}^{(s)}(P^{(s+1)},\Llambda,Y_{n-s+1})$ for $2\le j\le s$,
then computes $f_{s+1}$ modulo
$q_\lambda^{(s)}(P^{(s+1)},\Llambda,Y_{n-s+1})$, and finally
computes $g_{s+1}$. All these steps require
$O\big((\mathcal{T}+n)\mathcal{U}(\delta_s)\big)$ arithmetic
operations in ${\K}(\Llambda)$, storing at most
$O(\mathcal{S}\delta_s)$ elements of ${\K}(\Llambda)$. Each of
these arithmetic operations is performed in the power series ring
${\K}[\![\Llambda-\alpha]\!]$ at precision $d\delta_s+1$, and then
requires $O\big(\mathcal{U}(d\delta_s)\big)$ arithmetic operations
in $\K$, storing at most $O(d\delta_s)$ elements of $\K$.
Therefore, we conclude that the whole algorithm computing
$g_{s+1}$ requires $O\big((\mathcal{T}+n)\mathcal{U}(d\delta_s)
\mathcal{U}(\delta_s)\big)$ arithmetic operations in $\K$, storing
at most $O\big((\mathcal{S}+d)\delta_s^2\big)$ elements of $\K$.

Finally, the computation of
$q_{\Llambda}^{(s+1)}(P^{(s+1)},\Llambda)=g_{s+1}/\mathrm{gcd}
(g_{s+1},g_{s+1}')$ requires $O\big(\mathcal{U}(d\delta_s)\big)$
operations in $\K$, storing at most $O(d\delta_s)$ elements of
$\K$. This finishes the proof of the proposition.
\end{proof}
The algorithm underlying Proposition \ref{prop:minimal} is
essentially an extension to the finite field context of
\cite[Algorithm II.7]{GiLeSa01}. We have contributed further to
the latter by quantifying the probability of success of our
algorithm. Let us also remark that the complexity estimate of
Proposition \ref{prop:minimal} significantly improves that of
\cite[Proposition 1]{HeMaWa01}.
%
%
\subsection{Computing a geometric solution of \lfiberspone}
\label{subsec:shape_lemma}
In this section we exhibit an algorithm which computes a
parametrization of the variables $Y_{n-s+1}\klk Y_{n}$ by the
zeros of $q^{(s+1)}(P^{(s+1)},Y_{n-s})$, completing thus the
$s$--th recursive step of our main procedure computing a geometric
solution of the input variety $V$.

First we discuss how we obtain the parametrization of $Y_{n-s+1}$
by the zeros of $q^{(s+1)}(P^{(s+1)},Y_{n-s})$. Recall that such
parametrization is given by a polynomial $(\partial
q^{(s+1)}/\partial Y_{n-s})(P^{(s+1)},Y_{n-s}) Y_{n-s+1}-
v_{n-s+1}^{(s+1)}(P^{(s+1)}, Y_{n-s})\in {\K}[Y_{n-s},
Y_{n-s+1}]$, with $v_{n-s+1}^{(s+1)}(P^{(s+1)}, Y_{n-s})$ of
degree at most $\delta_{s+1}-1$.\spar

Let $\lambda_1,\lambda_2\in {\K}\setminus\{0\}$ satisfy the
conditions of Lemma \ref{lemma:moving_projection} and let
$\mathcal{L}_{\lambda_i}:=Y_{n-s}+\lambda_i Y_{n-s+1}$ for
$i=1,2$. Observe that the value $\lambda=0$ also satisfies the
conditions of Lemma \ref{lemma:moving_projection}. By Proposition
\ref{prop:minimal} we may assume that we have already computed the
minimal equations $q_{\Llambdaone}^{(s+1)}(P^{(s+1)},
\Llambdaone)$, $q_{\Llambdatwo}^{(s+1)}(P^{(s+1)},\Llambdatwo)$
and $q^{(s+1)}(P^{(s+1)},Y_{n-s})$ satisfied by
$\mathcal{L}_{\lambda_1}$, $\mathcal{L}_{\lambda_2}$ and $Y_{n-s}$
in $\cfq[\lfiberspone]$. Interpreting these polynomials as
elements of ${\K}[Y_{n-s},Y_{n-s+1}]$, assume further that
$\Llambdatwo$ separates the common zeros of
$q^{(s+1)}(P^{(s+1)},Y_{n-s})$ and $q_{\Llambdaone}^{(s+1)}
(P^{(s+1)},\Llambdaone)$. Arguing as above, we easily conclude that
there exists a nonzero polynomial $\widehat{E}_s\in\cfq[\Lambda]$
of degree at most $\delta^4$ such that, for any $\lambda_2$ with
$\widehat{E}_s(\lambda_2)\not=0$, the linear form $\Llambdatwo$
satisfies our last assumption.

In our subsequent argumentations we shall consider the following
(zero--dimen\-sional) $\K$--variety:
$$W_{\!s+\!1}\!:=\!\big\{(x_1,\!x_2)\!\in\!\A^2\!\!:\!q^{(s+\!1)}
(P^{(s+\!1)}\!,x_{1})\!=0,
q_{\mathcal{L}_{\lambda_i}}^{(s+1)}(P^{(s+\!1)}\!,
x_1+\lambda_ix_2)\!=0\mbox{ for }i\!=\!1,2\big\}.$$
Let $\widehat{\pi}_{s}^*:\lfiberspone\to\A^2$ the projection mapping
induced by $Y_{n-s},Y_{n-s+1}$. Observe that
$\widehat{\pi}_{s}^*(\lfiberspone)\subset W_{s+1}$ holds.
Furthermore, since $\Llambdatwo$ separates the common zeros of
$q^{(s+1)}(P^{(s+1)},Y_{n-s})$ and
$q_{\Llambdaone}^{(s+1)}(P^{(s+1)},\Llambdaone)$, and
$q_{\Llambdatwo}^{(s+1)}(P^{(s+1)},\Llambdatwo)$
vanishes in the set
$\Llambdatwo\big(\widehat{\pi}_{s}^*(\lfiberspone)\big)$ (of
cardinality $\delta_{s+1}$) and has degree $\delta_{s+1}$, we
conclude that $W_{s+1}=\widehat{\pi}_{s}^*(\lfiberspone)$ holds.

Our intention is to reduce the computation of
$v_{n-s+1}^{(s+1)}(P^{(s+1)},Y_{n-s})$ to $\gcd$ computations over
suitable field extensions of $\K$. From our previous argumentation
and the fact that $Y_{n-s}$ separates the points of $\lfiberspone$,
it follows that $Y_{n-s}$ also separates the points of $W_{s+1}$.
Then, applying the classical Shape Lemma to this
(zero--dimensional) $\K$--variety (see e.g. \cite{CoLiOS98}), we
see that there exists a polynomial $w_{n-s+1}\in {\K}[Y_{n-s}]$ of
degree at most $\delta_{s+1}-1$ such that
$Y_{n-s+1}-w_{n-s+1}(Y_{n-s})$ vanishes on the variety $W_{s+1}$.

Let $\alpha\in\cfq$ be an arbitrary root of
$q^{(s+1)}(P^{(s+1)},Y_{n-s})$ and let $\beta:=w_{n-s+1}(\alpha)$.
Then the fact that $Y_{n-s}$ separates the points of $W_{s+1}$
shows that $P:=(\alpha,\beta)$ is the only point of $W_{s+1}$ for
which $Y_{n-s}(P)=\alpha$ holds. Hence, $Y_{n-s+1}=\beta$ is the
only common root of
$q^{(s+1)}_{\mathcal{L}_{\lambda_1}}(P^{(s+1)},
\alpha+\lambda_1Y_{n-s+1})$ and
$q^{(s+1)}_{\mathcal{L}_{\lambda_2}}(P^{(s+1)},
\alpha+~\lambda_2Y_{n-s+1})$. Furthermore, the assumption on
$\lambda_2$ implies that
$q^{(s+1)}_{\mathcal{L}_{\lambda_2}}(P^{(s+1)},
\alpha+\lambda_2Y_{n-s+1})$ is squarefree. Therefore, we conclude
that the following identity holds in
${\K}(\alpha)[Y_{\!n-\!s+\!1}]$:
\begin{equation}
\label{eq:gcd}\gcd\!\Big(q^{(s+\!1)}_{\Llambdaone}(P^{(s+\!1)}\!,
\alpha\!+\!\lambda_1Y_{\!n-s+\!1}),
q^{(s+\!1)}_{\Llambdatwo}(P^{(s+\!1)}\!,\alpha\!+\!
\lambda_2Y_{\!n-s+\!1})\Big)\!\!=\!Y_{\!n-s+\!1}\!-\!\beta.
\end{equation}

Let $q^{(s+1)}(P^{(s+1)},Y_{n-s})= q_1\cdots q_h$ be the
irreducible factorization of the polynomial
$q^{(s+1)}(P^{(s+1)},Y_{n-s})$ in ${\K}[Y_{n-s}]$. Observe that
every irreducible factor $q_j$ represents an $\K$--irreducible
component $\mathcal{C}_j$ of $W_{s+1}$. Let $\alpha_j \in \cfq$ be
an arbitrary root of $q_j$. Taking into account the field
isomorphism ${\K}(\alpha_j)\simeq {\K}[Y_{n-s}]/
\big(q_j(Y_{n-s})\big)$, from identity (\ref{eq:gcd}) we conclude
that there exists $v_j\in {\K}[Y_{n-s}]$ of degree at most $\deg
q_j-1$ such that the following identity holds in
$\left({\K}[Y_{n-s}]/\big(q_j(Y_{n-s})\big)\right)[Y_{n-s+1}]$:
\begin{equation}
\label{eq:gcd2}
\gcd\!\!\Big(q^{(s+\!1)}_{\Llambdaone}(\!P^{(\!s+\!1)}\!,
Y_{\!n-\!s}\!+\!\lambda_1Y_{\!n-s+\!1}),
q^{(s+\!1)}_{\Llambdatwo}(\!P^{(\!s+\!1)}\!,Y_{\!n-\!s}\!+\!
\lambda_2Y_{\!n-s+\!1})\!\Big)\!\!=\!Y_{\!n-\!s+\!1}\!-\!
v_j(Y_{\!n-\!s}).\end{equation}
Let us fix $j\in\{1\klk h\}$. From the B\'ezout identity we deduce
that the congruence relation $Y_{n-s+1}-v_j(Y_{n-s})\equiv 0\mbox{
mod } I(\mathcal{C}_j)$ holds. This implies that $q_j'\cdot
(Y_{n-s+1}-v_j)$ belongs to the ideal $I(\mathcal{C}_j)$ for $1\le
j\le h$, and hence, that $q_j'\big(\prod_{i\not=j}
q_i\big)(Y_{n-s+1}-v_j)$ belongs to the ideal $I(W_{s+1})\subset
I(\lfiberspone)$ for $1\le j\le h$.

Let
\begin{equation}\label{eq:param_vnspone}
v_{n-s+1}^{(s+1)}(P^{(s+1)},Y_{n-s}):=\sum_{1\le j\le h}q_j'v_j
\big(\prod_{i\not=j} q_i\big)\ \mathrm{mod}\
q^{(s+1)}(P^{(s+1)},Y_{n-s}).
\end{equation}
By construction we have that
$v_{n-s+1}^{(s+1)}(P^{(s+1)},Y_{n-s})$ is an element of
${\K}[Y_{n-s}]$ of degree at most $\delta_{s+1}-1$. Furthermore,
our previous argumentation shows that $(\partial
q^{(\!s+\!1)}\!\!/\partial Y_{\!n-s+\!1})(\!P^{(\!s+\!1)}\!,
Y_{\!n-\!s})Y_{\!n-\!s+\!1}- v_{\!n-s+\!1}^{(\!s+\!1)}
(\!P^{(\!s+\!1)}\!,Y_{\!n-\!s})\!\!= \!\!\sum_{j=1}^h
\!q_{\!j}'\big(\!\prod_{i\not=j}\! q_i\big)(Y_{\!n-s+\!1}-v_j)$
belongs to the ideal $I(\lfiberspone)$, and hence it represents the
parametrization of $Y_{n-s+1}$ by the zeros of
$q^{(s+1)}(P^{(s+1)},Y_{n-s})$ we are looking for. \spar

Now we discuss how we can obtain the polynomial
$v_{n-s+j}^{(s+1)}(P^{(s+1)},Y_{n-s})$, which parametrizes
$Y_{n-s+j}$ by the zeros of $q^{(s+1)}(P^{(s+1)},Y_{n-s})$,
for $2\le j\le s$.
For this purpose, we assume that we are given the
polynomials $v_{n-s+j}^{(s)}(P^{(s+1)},Y_{n-s},Y_{n-s+1})$ $(2\le
j\le s)$ which form the geometric solution of the lifting curve
$\lcurvespone$.

Let $(\partial q^{(s+1)}/\partial
Y_{n-s+1})^{-1}(P^{(s+1)},Y_{n-s})\in {\K}[Y_{n-s}]$ denote the
inverse of the polynomial $(\partial q^{(s+1)}/\partial
Y_{n-s+1})(P^{(s+1)},Y_{n-s})$ modulo
$q^{(s+1)}(P^{(s+1)},Y_{n-s})$, and let
$w_{n-s+1}^{(s+1)}(P^{(s+1)},Y_{n-s})\!:=(\partial
q^{(s+1)}/\partial Y_{n-s+1})^{-1}(P^{(s+1)},Y_{n-s})\,
v_{n-s+1}^{(s+1)}(P^{(s+1)},Y_{n-s})$. Observe that
$Y_{n-s+1}-w_{n-s+1}^{(s+1)}(P^{(s+1)},Y_{n-s})$ belongs to the
ideal $I(\lfiberspone)$. Our intention is to use this
parametrization to ``clean" the variable $Y_{n-s+1}$ out of the
polynomials $v_{n-s+j}^{(s+1)}(P^{(s+1)},Y_{n-s},Y_{n-s+1})$.

For this purpose, we observe that $q^{(s)}\big(P^{(s+1)},Y_{n-s},
w_{n-s+1}^{(s+1)}(P^{(s+1)},Y_{n-s})\big)$ and $(\partial
q^{(s)}\!\!/\partial Y_{n-s+1})\big(P^{(s+1)},Y_{n-s},
w_{n-s+1}^{(s+1)}(P^{(s+1)},Y_{n-s})\big)Y_{n-s+j} -
v_{n-s+j}^{(s+1)}\big(P^{(s+1)}\!,Y_{n-s},$
$w_{n-s+1}^{(s+1)}(P^{(s+1)}, Y_{n-s})\big)$ $(2\le j\le s)$
belong to the ideal $I(\lfiberspone)$. Furthermore, we have that
the polynomial $(\partial q^{(s)}\!\!/\partial
Y_{\!n-s+\!1})\big(P^{(s+\!1)}\!,Y_{\!n-\!s},
w_{n-\!s+\!1}^{(s+\!1)}(P^{(s+\!1)}\!,Y_{\!n-\!s})\big)$ is not a
zero divisor of
${\K}[Y_{n-s}]/\big(q^{(s+1)}(P^{(s+1)},Y_{n-s})\big)$, because
otherwise the discriminant $\rho^{(s)}(P^{(s+1)},Y_{n-s})$ would
have common roots with $q^{(s+1)}(P^{(s+1)},Y_{n-s})$,
contradicting thus condition $(iii)$ of Theorem
\ref{theorem:simNoether}. Therefore, its inverse $h_{n-s+1}$
modulo $q^{(s+1)}(P^{(s+1)},Y_{n-s})$ is well--defined element of
${\K}[Y_{n-s}]$, and the polynomial $Y_{n-s+j} - h_{n-s+1}\cdot
v_{n-s+j}^{(s+1)}\big(P^{(s+1)},Y_{n-s},
w_{n-s+1}^{(s+1)}(P^{(s+1)}, Y_{n-s})\big)$ belongs to
$I(\lfiberspone)$ for $2\le j\le s$. Therefore, if we let
\begin{equation}\label{eq:def_w}w_{n-s+j}\!:=\!h_{n-s+1}\cdot
v_{n-s+j}^{(s+1)}\big(P^{(s+\!1)}\!,Y_{n-s},
w_{n-s+1}^{(s+1)}(P^{(s+\!1)}\!,Y_{n-s})\big)\  (2\le j\le s),
\end{equation}
we see that $Y_{n-s+j}-w_{n-s+j}$ belongs to $I(\lfiberspone)$ for
$2\le j\le s$. Multiplying $w_{n-s+j}$ by $(\partial
q^{(s+1)}/\partial Y_{n-s+1})(P^{(s+1)},Y_{n-s})$ for $2\le j\le
s$, and reducing modulo $q^{(s+1)}(P^{(s+1)},Y_{n-s})$, we obtain
the polynomials $v_{n-s+j}^{(s+1)}\in {\K}[Y_{n-s}]$ $(2\le j\le
s)$ we are looking for.\spar

In our next result we exhibit an algorithm computing the
polynomials \linebreak $q^{(s+1)}(P^{(s+1)},Y_{n-s}),
v_{n-s+j}^{(s+1)}(P^{(s+1)},Y_{n-s})\in {\K}[Y_{n-s}]$ $(1\le j\le
s)$, which form a geometric solution of $\lfiberspone$.
\begin{proposition}\label{prop:recursive_geo_sol}
There exists a probabilistic Turing machine $M$ which has as input
\begin{itemize}
\item a \slp\ using space $\mathcal{S}$ and time $\mathcal{T}$
which represents the polynomial $F_{s+1}$,
\item the polynomials $q^{(s)}(P^{(s+\!1)}\!,Y_{n-s},Y_{n-s+1})$ and
$v_{n-s+j}^{(s)}(P^{(s+\!1)}\!,Y_{n-s},Y_{n-s+1})$ ($2\le j\le
s)$, which form the geometric solution of the lifting curve
$\lcurvespone$ computed in Proposition
\ref{prop:lifting_to_curve},
\end{itemize}
and outputs a geometric solution of the lifting fiber
$\lfiberspone$. The Turing machine $M$ runs in space
$O\big((\mathcal{S}+n+d)\delta^2\log (q\delta)\big)$ and time
$O\big((\mathcal{T}+n)\mathcal{U}(\delta)\big(\mathcal{U}(d\delta)
+\log (q\delta)\big)\mathcal{U}(\log (q\delta))\big)$, and outputs
the right result with probability at least $1-1/60n$.
\end{proposition}
\begin{proof}
Let $E_s$ be the polynomial of the statement of Lemma
\ref{lemma:moving_projection} and let $\widehat{E}_s$ be the
polynomial introduced at the beginning of this section. Recall
that $\deg E_s\le 4\delta^3$ and $\deg \widehat{E}_s\le \delta^4$
hold. Let $\lambda_1,\lambda_2$ two distinct values of $\K$
randomly chosen and let
$\mathcal{L}_{\lambda_i}:=Y_{n-s}+\lambda_i Y_{n-s+1}$ ($i=1,2$).
Applying Theorem \ref{th:Zippel_Schwartz} we conclude that
$E_s(\lambda_1) E_s(\lambda_2)\widehat{E}_s(\lambda_2)\not=0$
holds with probability at least $1-1/72n^3$. Suppose that this is
the case. Then, applying the algorithm underlying Proposition
\ref{prop:minimal}, we conclude that the minimal equations
$q^{(s+1)}(P^{(s+1)},Y_{n-s}),
q^{(s+1)}(P^{(s+1)},\mathcal{L}_{\lambda_i})$ $(i=1,2)$ satisfied
by $Y_{n-s},\mathcal{L}_{\lambda_i}$ ($i=1,2$) in
${\K}[\lfiberspone]$ can be computed by a probabilistic Turing
machine which runs in space $O\big((\mathcal{S}+d)\delta_s^2\log
(q\delta)\big)$ and time
$O\big((\mathcal{T}+n)\mathcal{U}(d\delta_s)
\mathcal{U}(\delta_s)\mathcal{U}(\log (q\delta))\big)$, with
probability of success at least $1-1/15n^3$.

Next we compute the irreducible factorization
$q^{(s+1)}(P^{(s+1)},Y_{n-s})=q_1\cdots q_h$ of
$q^{(s+1)}(P^{(s+1)},Y_{n-s})$ in ${\K}[Y_{n-s}]$. From
\cite[Corollary 14.30]{GaGe99} we conclude that such factorization
can be computed with space $O(\delta_{s+1}^2\!\log(q\delta))$ and
time $O\big(\!\log(n)\big(\mathcal{U}(\delta_{\!s+\!1}^2)$ $
+\mathcal{U}( \delta_{\!s+\!1})\log (q\delta)\big)\mathcal{U}(\log
(q\delta))\big)$, with probability of sucess at least $1-1/16n^3$.

Then we compute the polynomials $v_1\klk v_h$ of (\ref{eq:gcd2})
and the polynomial $v_{n-s+1}^{(s+1)}$ of (\ref{eq:param_vnspone}),
using the EEA (see e.g. \cite{BiPa94}, \cite{GaGe99}). According
to \cite[Corollary 11.16]{GaGe99}, this step can be done
deterministically using space $O(\delta_s\delta_{s+1}\log
(q\delta))$ and time \linebreak
$O\big(\delta_{s+1}\mathcal{U}(\delta_s)\mathcal{U}(\log
(q\delta))\big)$. Finally, we compute the polynomials $h_{n-s+1}$
and $w_{n-s+j}$ $(2\le j\le s)$ of (\ref{eq:def_w}), and the
polynomials $v_{n-s+j}^{(s+1)}(P^{(s+1)},Y_{n-s})$ for $2\le j\le
s$, with the same asymptotic complexity estimates. Adding the
complexity and probability estimates of each step, we easily deduce
the statement of the proposition.
\end{proof}
The algorithm underlying Proposition \ref{prop:recursive_geo_sol}
extends to the positive characteristic case the algorithms of
\cite{HeMaWa01} and \cite{GiLeSa01}, having a better asymptotic
complexity estimate (in terms of the number of arithmetic
operations performed) than \cite{HeMaWa01}, and a similar
complexity estimate as \cite{GiLeSa01}. We also contribute to the
latter by providing estimates on the probability of success of the
algorithm, which are not present in \cite{GiLeSa01}. Finally, let
us also remark that by means of our preprocessing we have
significantly simplified both the algorithms of \cite{HeMaWa01}
and \cite{GiLeSa01}.
%
%
\subsection{A $\K$--definable geometric solution of $V$}
\label{subsec:geo_sol_final}
Now we have all the ingredients necessary to describe our
algorithm computing the $\K$--definable geometric solution of our
input variety $V:=V_r$. We recall that $\K$ is a field extension of
$\fq$ of cardinality greater than $60n^4d\delta^4$. Let
$(\lambda,\gamma,P)$ be a point randomly chosen in the set
${\K}^{n^2}\times {\K}^n\times {\K}^{n-r}$. Theorem
\ref{th:Zippel_Schwartz} shows that $B(\lambda,\gamma,P)$ does not
vanish with probability at least $14/15$, where $B$ is the
polynomial defined at the beginning of Section
\ref{section:computation_geo_sol}. Assume that we have chosen such
a point and let $(Y_1\klk Y_n):=\lambda X+\gamma$ and $P:=(p_1\klk
p_{n-r})$. Then $Y_1\klk Y_n$ and $P^{(s)}:=(p_{1}\klk p_{n-s})$
satisfy the conditions of Theorem \ref{theorem:simNoether} for
$1\le s\le r-1$.

Therefore, we may recursively apply, for $1\le s\le r-1$, the
algorithms underlying Propositions \ref{prop:lifting_to_curve} and
\ref{prop:recursive_geo_sol}, which compute a geometric solution of
the lifting curve $\lcurvespone$ and of the lifting fiber
$\lfiberspone$ respectively. In this way, at the end of the
$(r-1)$--th recursive step we obtain a geometric solution of the
lifting fiber $\lfiberr$. Taking into account the complexity and
probability estimates of Propositions \ref{prop:lifting_to_curve}
and \ref{prop:recursive_geo_sol}, we easily deduce the following
result:
\begin{theorem}\label{th:geo_sol} There exists a
probabilistic Turing machine $M$, which takes as input a \slp\
which represents the input polynomials $\fo{r}$ with space
$\mathcal{S}$ and time $\mathcal{T}$, and outputs a geometric
solution of the lifting fiber $\lfiberr$. The Turing machine $M$
runs in space $O\big((\mathcal{S}+n+d)\delta^2\log (q\delta)\big)$
and time $O\big((n\mathcal{T}+n^5)\mathcal{U}(\delta)
\big(\mathcal{U}(d\delta)+\log (q\delta)\big)\mathcal{U}(\log
(q\delta))\big)$ and outputs the right result with probability at
least $1-1/12$.
\end{theorem}
The complexity estimate of Theorem \ref{th:geo_sol} significantly
improves the $O(d^{n^2})$ complexity estimate of \cite{HuWo99},
the $O(d^{2r})$ estimate of \cite{HuWo00}, and the estimates of
the algorithms of the so--called Gr\"obner solving. Further, let
us remark that, combining the algorithm underlying Theorem
\ref{th:geo_sol} with techniques of $p$--adic lifting, as those of
\cite{GiLeSa01}, for a ``lucky" choice of prime number $p$, one
obtains an efficient probabilistic algorithm computing the
geometric solution of an equidimensional variety over $\Q$ given
by a reduced regular sequence.
%
%
\section{An $\fq$--definable lifting fiber of $V$}
\label{subsec:fq_def_geo_sol}
Let notations and assumptions be as Section
\ref{subsec:geo_sol_final}. In this section we obtain a geometric
solution of an $\fq$--definable lifting fiber of $V$. For this
purpose, we shall homotopically deform the $\K$--definable
geometric solution of the lifting fiber
$\lfiberr:=\pi_r^{-1}(P^{(r)})$, computed in the previous section,
into a geometric solution of an $\fq$--definable lifting fiber
$\pi^{-1}(Q)$ of the linear projection mapping $\pi:V\to \A^{n-r}$
determined by suitable linear forms $Z_1\klk Z_{n-r}\in\fq[\xon]$.

Let  $\Lambda $ be an $(n-r+1)\times n$ matrix of indeterminates.
For $1\le i\le n-r+1$, let
$\Lambda^{(i)}:=(\Lambda_{i1},\dots,\Lambda_{in})$ denote its
$i$--th row and let $\Lambda^{(1\klk i)}$ denote the $i\times n$
submatrix of $\Lambda$ consisting of the first $i$ rows of $\Lambda$.
Let $\Gamma:=(\Gamma_1\klk \Gamma_{n})$ be a vector of
indeterminates, and let $\widetilde{Y}:= (\widetilde{Y}_1\klk
\widetilde{Y}_{n-r+1}):=\Lambda X + \Gamma$.

Let $\widehat{B}\in\cfq[\Lambda,\Gamma,\widetilde{Y}_1\klk
\widetilde{Y}_{n-r}]$ be the polynomial of Corollary
\ref{coro:Noether_V_r}, and let $B':=\det(\Delta_1)\det(\Delta_2)
\widehat{B}$, where $\Delta_1$ is the $n\times n$ matrix
which has $\Lambda^{(1\dots n-r)}$ as its upper $(n-r)\times n$
submatrix, and the coefficients of the linear forms $Y_{n-r+1}\klk
Y_n$ in its last $r$ rows, and $\Delta_2$ is the $n\times n$
matrix having $\Lambda^{(1\dots n-r+1)}$ as its upper $(n-r+1)\times n$
submatrix, and the coefficients of $Y_{n-r+2}\klk Y_n$ in its last
$r-1$ rows. Observe
that $\deg B'\le 2(n-r+2)nd\delta_r^2$ holds.

Suppose that $q>8n^2d\delta_r^4$ holds, and let be given a point
$(\nu,\eta,Q)\in\fq^{(n-r+1)n}\times \fq^{n-r+1}\times\fq^{n-r}$
such that $B'(\nu,\eta,Q)\not=0$. Theorem \ref{th:Zippel_Schwartz}
shows that such a point $(\nu,\eta,Q)$ can be randomly chosen
in the set $\fq^{(n-r+1)n}\times \fq^{n-r+1}\times\fq^{n-r}$
with probability of success at least $1-1/16$.

Let $\nu:=\nu^{(1\klk n-r+1)}$, $\eta:=(\eta_1\klk \eta_{n-r+1})$,
$Q:=(q_1\klk q_{n-r})$ and $Z:=(Z_1\klk Z_{n-r+1}):=\nu X+\eta$.
The condition $\det(\Delta_1\cdot\Delta_2)(\nu)\not=0$ implies that
the sets of linear forms $Z_1\klk Z_{n-r},Y_{n-r+1}\klk Y_n$ and
$Z_1\klk Z_{n-r+1},Y_{n-r+2}\klk Y_n$ induce linear changes of
coordinates. Furthermore, from the condition
$\widehat{B}(\nu,\eta,Q)\not=0$ and Corollary
\ref{coro:Noether_V_r}, we conclude that the linear projection mapping
$\pi:V\to\A^{n-r}$ defined by $Z_1\klk Z_{n-r}$ is a finite
morphism, $Q\in \fq^{n-r}$ is a lifting point of $\pi$ and
$Z_{n-r+1}$ is a primitive element of the lifting fiber
$V_Q:=\pi^{-1}(Q)$.

Let $(\lambda, \gamma, P)\in \K^{n^2}\times \K^n \times \K^{n-r}$
be the point fixed at the beginning of Section \ref{section:computation_geo_sol},
which yields the linear forms $Y:=(Y_1\klk Y_n):=\lambda X+\gamma$ and the point
$P^{(r)}\in\K^{n-r}$. Write
$\gamma:=(\gamma_1\klk \gamma_n)$ and $P:=(p_1\klk p_{n-r})$

Let  T be a new indeterminate, and let $\widehat{\Lambda} \in
\cfq(T)^{n\times n}$ and $\widehat{\Gamma}\in \cfq(T)^{n}$ be
the matrix and column vector defined in the following way:
$$\begin{array}{rcl}
\widehat{\Lambda}&:=&(1-T)\lambda\,\,+T\Delta_1(\nu^{(1\dots n-r)}),\\
\widehat{\Gamma}&:=& (1-T)\gamma^t+T (\eta_1\klk\eta_{n-r},\gamma_{n-r+1}
\klk \gamma_n)^t,\end{array}$$
where $\nu^{(1 \klk n-r)}$ denotes the $(n-r)\times n$ matrix
consisting of the first $n-r$ rows of $\nu$ and the symbol $^t$
denotes transposition. Let $\widehat{\Lambda}^{(1 \klk n-r)}$
and $\widehat{\Gamma}^{(1\klk n-r)}$ denote the $(n-r)\times n$
submatrix of $\widehat{\Lambda}$ consisting of the first $n-r$
of $\widehat{\Lambda}$ and the $(n-r)$--dimensional vector
consisting of the first $n-r$ entries of $\widehat{\Gamma}$
respectively.

Let $W$ be the subvariety of $\overline{\fq(T)}^{n}$ defined by
the set of common zeros of $F_1\klk F_r$. Let
$\widehat{Z}:=(\widehat{Z}_1\klk \widehat{Z}_n)
:=\widehat{\Lambda} X + \widehat{\Gamma}$ and
$\widehat{P}:=(\widehat{p}_1\klk
\widehat{p}_{n-r}):=(1-T)P+ TQ$. Since $\widehat{\Lambda}$ is an
invertible element of $\cfq(T)^{n\times n}$, we have that
$X=\widehat{\Lambda}^{-1} (\widehat{Z}-\widehat{\Gamma})$
holds, and hence $\widehat{F}_j(T,\widehat{Z}):=
F_j(\widehat{\Lambda}^{-1}(\widehat{Z}-\widehat{\Gamma}))$
is a well--defined element of $\cfq(T)[\widehat{Z}_1\klk
\widehat{Z}_n]$ for $1\le j\le r$. Observe that the point
$(\widehat{\Lambda},\widehat{\Gamma},\widehat{P})\in
\overline{\fq(T)}^{n^2}\times\overline{\fq(T)}^{n}\times
\overline{\fq(T)}^{n-r}$ does not annihilates the polynomial
$\widehat{B}$ of the statement of Corollary
\ref{coro:Noether_V_r}. Therefore, applying Corollary
\ref{coro:Noether_V_r}, replacing the field $\cfq$ by
$\overline{\fq(T)}$, we conclude that
$\overline{\fq(T)}[\widehat{Z}_1\klk
\widehat{Z}_{n-r}]\hookrightarrow \overline{\fq(T)}[X]/(F_1\klk
F_r)$ is an integral ring extension, $\widehat{P}$ is a lifting
point of the linear projection mapping $\pi^e:W\to
\overline{\fq(T)}^{n-r}$ defined by $\widehat{Z}_1\klk
\widehat{Z}_{n-r}$, and $\widehat{Z}_{n-r+1}=Y_{n-r+1}$ is a
primitive element of the (zero--dimensional) lifting fiber
$W_{\widehat{P}}:=(\pi^e)^{-1}(\widehat{P})$.

Let $\widehat{q}_{\widehat{Z}_{n-r+1}}(\widehat{P},\widehat{Z}_{n-r+1})
\in\overline{\fq(T)}[\widehat{Z}_{n-r+1}]$ denote the minimal equation
satisfied by
$\widehat{Z}_{n-r+1}$ in $\overline{\fq(T)}[W_{\widehat{P}}]$.
By the $\K(T)$--definability of $W_{\widehat{P}}$ and
$\widehat{Z}_{n-r+1}$, we see that $\widehat{q}_{
\widehat{Z}_{n-r+1}}(\widehat{P},\widehat{Z}_{n-r+1})$ belongs to
$\K(T)[\widehat{Z}_{n-r+1}]$.
Furthermore, our choice of $\widehat{P}$ and
$\widehat{Z}_1\klk \widehat{Z}_{n-r+1}$ implies that
$\widehat{q}_{ \widehat{Z}_{n-r+1}}(\widehat{P},\widehat{Z}_{n-r+1})$ is a
separable element of $\K(T)[\widehat{Z}_{n-r+1}]$ of degree $\delta_r$. Let
$\widehat{\rho}\in\K[T]$ be the product of the denominator of
$\widehat{q}_{\widehat{Z}_{n-r+1}}(\widehat{P},\widehat{Z}_{n-r+1})$ and the
numerator of its discriminant. In order to perform the homotopic
deformation mentioned at the beginning of this section, we need the following
preliminary result:
\begin{lemma}\label{lemma:cond_fq_def}
Let $\widehat{I}_{\widehat{P}}$ be the ideal of
${\K}[T]_{\widehat{\rho}}[Y_{\!n-r+1}\klk\!Y_{\!n}]$
generated by the polynomials
$\widehat{F}_1(\widehat{P},Y_{n-r+1}\klk Y_{n})\klk
\widehat{F}_r(\widehat{P},Y_{n-r+1}\klk Y_{n})$. Then
\begin{equation}\label{eq:integral_ext_fq_def}
{\K}[T]_{\widehat{\rho}}\hookrightarrow
{\K}[T]_{\widehat{\rho}}[Y_{\!n-r+1}\klk\!Y_{\!n}]/
\widehat{I}_{\widehat{P}}\end{equation}
is an integral ring extension, the polynomials
$\widehat{F}_j(\widehat{P},Y_{n-r+1}\klk Y_{n})$
$(1\le j\le r)$ form a regular
sequence of ${\K}[T]_{\widehat{\rho}}[Y_{n-r+1}\klk Y_{n}]$ and
generate a radical ideal $\widehat{I}_{\widehat{P}}$ of
${\K}[T]_{\widehat{\rho}}[Y_{n-r+1}\klk Y_{n}]$, and the ring
${\K}[T]_{\widehat{\rho}}[Y_{n-r+1}\klk Y_{n}]/
\widehat{I}_{\widehat{P}}$
is a free $\K[T]_{\widehat{\rho}}$--module of degree $\delta_r$.
\end{lemma}
\begin{proof}
By the remarks above, we see that $\widehat{q}_{\widehat{Z}_{n-r+1}}
(\widehat{P},\widehat{Z}_{n-r+1})\in\K[T]_{\widehat{\rho}}
[\widehat{Z}_{n-r+1}]$ is an integral dependence equation for the
coordinate function $\widehat{z}_{n-r+1}$ induced by $\widehat{Z}_{n-r+1}$
in the ring extension (\ref{eq:integral_ext_fq_def}). We conclude
that $\K[T]_{\widehat{\rho}}\hookrightarrow\K[T]_{\widehat{\rho}}
[\widehat{z}_{n-r+1}]$ is an integral ring extension.

Let $\xi_1\klk \xi_n$ denote the coordinate functions of
$\K[T]_{\widehat{\rho}}[Y_{n-r+1}\klk Y_n]/\widehat{I}_{\widehat{P}}$
induced by $X_1\klk X_n$. Arguing as in eq. (\ref{equation:chow3})
of the proof of Proposition \ref{prop:chow}, we conclude that there
exists polynomials $\widehat{P}_1\klk \widehat{P}_n\in\K[T]_{\widehat{\rho}}
[\widehat{Z}_{n-r+1}]$ such that $\xi_k=\widehat{P}_k(\widehat{z}_{n-r+1})$
holds for $1\le k\le n$. This shows that $\K[T]_{\widehat{\rho}}
[\widehat{z}_{n-r+1}]\to\K[T]_{\widehat{\rho}}[\xi_1\klk \xi_n]$ $=
\K[T]_{\widehat{\rho}}[Y_{n-r+1}\klk Y_n]/\widehat{I}_{\widehat{P}}$
is an integral ring extension and, combined with fact that
$\K[T]_{\widehat{\rho}}\hookrightarrow\K[T]_{\widehat{\rho}}
[\widehat{z}_{n-r+1}]$ is an integral ring extension, proves our
first assertion.

Now we show that
$\widehat{F}_1(\widehat{P},Y_{n-r+1}\klk Y_{n})\klk
\widehat{F}_r(\widehat{P},Y_{n-r+1}\klk Y_{n})$ form a regular
sequence of ${\K}[T]_{\widehat{\rho}}[Y_{n-r+1}\klk Y_{n}]$ and
the ideal $\widehat{I}_{\widehat{P}}$ is radical. Arguing by
contradiction, suppose that there exists $1\le j\le r$ such that
$\widehat{F}_j(\widehat{P},Y_{n-r+1}\klk Y_{n})$ is a zero
divisor modulo $\widehat{F}_1(\widehat{P},Y_{n-r+1}\klk Y_{n})\klk
\widehat{F}_{j-1}(\widehat{P},Y_{n-r+1}\klk Y_{n})$. By specialization
in $T=0$ we conclude that $\widehat{F}_j(P^{(r)},Y_{n-r+1}\klk Y_{n})$
is a zero divisor modulo
$\widehat{F}_1(P^{(r)},Y_{n-r+1}\klk Y_{n})\klk
\widehat{F}_{j-1}(P^{(r)},Y_{n-r+1}\klk Y_{n})$, contradicting thus
Lemma \ref{lemma:sucesionregular}. Furthermore, a similar argument
shows that the Jacobian determinant
$\det\big(\partial\widehat{F}_i(\widehat{P},Y_{n-r+1}\klk Y_{n})/
\partial Y_{n-r+j}\big)_{1\le i,j\le r}$ is not a zero divisor
modulo $\widehat{F}_1(\widehat{P},Y_{n-r+1}\klk Y_{n})\klk
\widehat{F}_{r}(\widehat{P},Y_{n-r+1}\klk Y_{n})$. Hence,
from \cite[Theorem 18.15]{Eisenbud95} we deduce that the ideal
$\widehat{I}_{\widehat{P}}$ is radical.

Our previous assertions imply that
${\K}[T]_{\widehat{\rho}}[Y_{n-r+1}\klk Y_{n}]/\widehat{I}_{\widehat{P}}$
is a free ${\K}[T]_{\widehat{\rho}}$--module or rank at most $\delta_r$.
Since $\widehat{q}_{\widehat{Z}_{n-r+1}}(\widehat{P},\widehat{Z}_{n-r+1})$
is the minimal dependence equation satisfied by $\widehat{z}_{n-r+1}$
in the extension (\ref{eq:integral_ext_fq_def}), we conclude that
the rank of ${\K}[T]_{\widehat{\rho}}[Y_{n-r+1}\klk Y_{n}]/\widehat{I}_{\widehat{P}}$
as ${\K}[T]_{\widehat{\rho}}$ is exactly $\delta_r$. This finishes the
proof of the lemma.
\end{proof}

Let $\widehat{V}\subset\A^{r+1}$ be the affine equidimensional
variety defined by the set of common zeros of $\widehat{F}_1(\widehat{P},Y_{n-r+1}\klk
Y_{n})\klk$ $\widehat{F}_r(\widehat{P},Y_{n-r+1}\klk Y_{n})$
and let $\widehat{\pi}:\widehat{V}\to\A^1$ be the mapping
induced by the projection onto the coordinate $T$. Lemma
\ref{lemma:cond_fq_def} implies that $\widehat{V}$ has dimension
1 and degree $\delta_r$, and $\widehat{\pi}$ is a finite morphism.
Furthermore, taking into
account the equalities $\widehat{V}\cap \{T=0\}=\{0\} \times \lfiberr$ and
$\widehat{V}\cap \{T=1\}=\{1\} \times V_Q$, we conclude that $T=0,1$
are lifting points of the morphism $\widehat{\pi}$.
Therefore, applying a suitable variant of the Newton--Hensel
procedure of Section \ref{subsec:lifting_to_curve}, we obtain
a geometric solution of the lifting fiber $V_Q$. This is the
content of our next result:
\begin{proposition}
\label{prop:fq_geo_sol} Suppose that $q>8n^2d\delta_r^4$ holds.
Then there exists a probabilistic Turing machine $M$ which has as
input
\begin{itemize}
\item a \slp\ using space $\mathcal{S}$ and time $\mathcal{T}$
which represents the input polynomials $\fo{r}$,
\item the polynomials $q^{(r)}(P^{(r)},Y_{n-r+1})$,
$v_{n-s+j}^{(r)}(P^{(r)},Y_{n-r+1})$ ($2\le j\le s)$, which form
the geometric solution of the lifting fiber $\lfiberr$ computed in
Theorem \ref{th:geo_sol},
\end{itemize}
and outputs polynomials $q(Q,Z_{n-r+1})$,
$v_{n-s+j}(Q,Z_{n-r+1})$ ($2\le j\le s)$ which form a geometric
solution of the lifting fiber $V_Q$. The Turing machine $M$ runs
in space $O\big((\mathcal{S}+n)\delta_r^2\log(q\delta)\big)$ and
time $O\big((n\mathcal{T}+n^5)\mathcal{U}(\delta_r)^2
\mathcal{U}(\log (q\delta))\big)$ and outputs the right result
with probability at least 1-1/16.
\end{proposition}
\begin{proof}
Let $(\nu,\eta,Q)$ be a point randomly chosen in the set
$\fq^{(n-r+1)n}\times \fq^{n-r+1}\times\fq^{n-r}$. Let
$B'\in\cfq[\Lambda,\Gamma,\widetilde{Y}_1\klk \widetilde{Y}_{n-r}]$
be the polynomial introduced at the beginning of this section. Since
$\deg B'\le 2(n-r+2)nd\delta_r^2$ holds, from Theorem
\ref{th:Zippel_Schwartz} we conclude that $B'(\nu,\eta,Q)\not=0$
holds with probability at least $1-1/16$.

By the remarks before the statement of the proposition, we see that
$T=0,1$ are lifting points of the morphism
$\widehat{\pi}$. Then, applying the Newton--Hensel procedure of
\cite{Schost03}, we see that there exists a computation tree in
$\K$ computing polynomials $\widehat{q}(T,Y_{n-r+1})$,
$\widehat{v}_{n-r+j}(T,Y_{n-r+1})$ $(2\le j\le r)$ which form a
geometric solution of $\widehat{V}$. This computation tree
requires $O\big((n\mathcal{T}+n^5)\mathcal{U}(\delta_r)^2\big)$
operations in $\K$, using at most $O\big((\mathcal{S}+n)
\delta_r^2\big)$ arithmetic registers. Specializing these
polynomials into $T=1$ yields polynomials
$\widehat{q}(1,Y_{n-r+1})$, $\widehat{v}_{n-s+j}(1,Y_{n-r+1})$
$(2\le j\le r)$, which form a geometric solution of the lifting
fiber $\widehat{V}\cap \{T=1\}=\{1\}\times V_Q$ (and therefore
of $V_Q$), using $Y_{n-r+1}$ as primitive element.

Our next purpose is to compute a geometric solution of $V_Q$, using
$Z_{n-r+1}$ as primitive element. In order to to this, let
$\widehat{w}_{n-r+j}(1,Y_{n-r+1})\in\K[S]$ denote the remainder of
the product $(\partial\widehat{q}/\partial Y_{n-r+1})
(1,Y_{n-r+1})\cdot\widehat{v}_{n-r+j}(1,Y_{n-r+1})$ modulo
$\widehat{q}(1,Y_{n-r+1})$ for $2\le j\le r$. Observe that
$Y_{n-r+j}=\widehat{w}_{n-r+j}(1,Y_{n-r+1})$ holds in $\K[V_Q]$
for $2\le j\le r$. Write $Z_{n-r+1}=\alpha_1
Z_1+\cdots+\alpha_{n-r}Z_{n-r}+\alpha_{n-r+1}Y_{n-r+1}\plp\alpha_n
Y_n$. Then, it is easy to see that the minimal equation satisfied
by the linear form $Z_{n-r+1}+TY_{n-r+1}$ in
$\cfq[T]\otimes\cfq[V_Q]$ is given by
\begin{equation}
\label{eq:res_homot_projection} q_{Z_{n-r+1}+TY_{n-r+1}}(Q,S)=
\mathrm{Res}_U\Big(\widehat{q}(1,U),U-
\!\!\sum_{j=1}^{n-r}\!\alpha_jq_j+\alpha_{n-r+1}S+\!\!\!\!\!\!\!\!
\sum_{j=n-r+2}^n \!\!\!\!\alpha_j \widehat{w}_j(1,S)\Big).
\end{equation}
Following e.g. \cite{AlBeRoWo96} or \cite{Rouillier97} as in the proof of
Proposition \ref{prop:chow}, we
have the congruence relation
$q_{Z_{n-r+1}+TY_{n-r+1}}(Q,Z_{n-r+1})\equiv
q(Q,Z_{n-r+1})+\big(\partial q/\partial Z_{n-r+1}(Q,Z_{n-r+1})-
v_{n-r+1}(Q,Z_{n-r+1})\big)T$ modulo $(T^2)$, where
$q(Q,Z_{n-r+1})$ is the minimal equation of $Z_{n-r+1}$
in $\K[V_Q]$ and $(\partial q/\partial
Z_{n-r+1})(Q,Z_{n-r+1})Y_{n-r+1}=v_{n-r+1}(Q,Z_{n-r+1})$ holds in
$\K[V_Q]$.

We compute the right--hand--side term of
(\ref{eq:res_homot_projection}), up to order $T^2$, by interpolation
in the variable $S$, reducing thus the computation to $\delta_r$ resultants
of univariate polynomials of $\K[T]$ of degree at most $1$. Using
fast algorithms for univariate resultants and
interpolation over $\K$ (see e.g. \cite{BiPa94}, \cite{GaGe99}),
we conclude that the dense representation of $q(Q,S)$ and
$v_{n-r+1}(Q,S)$ can be deterministically computed with
$O(\delta_r\mathcal{U}(\delta_r))$ arithmetic operations over
$\K$, using at most $O(\delta_r^2)$ arithmetic registers.

Finally, there remains to compute the polynomials
$v_{n-r+j}(Q,Z_{n-r+1})\ (2\le j\le r)$ which parametrize
$Y_{n-r+j}$ by the zeros of  terms of $q(Q,Z_{n-r+1})$, i.e. such
that $(\partial q/\partial Z_{n-r+1})(Q,Z_{n-r+1})Y_{n-r+j}\equiv
v_{n-r+j}(Q, Z_{n-r+1})$ holds in $\K[V_Q]$. For this purpose, we
shall compute polynomials $w_{n-r+j}(Q,Z_{n-r+1})\ (1\le j\le r)$
of degree at most $\delta_r-1$ such that $Y_{n-r+j}\equiv
w_{n-r+j}(Q, Z_{n-r+1})$ holds in $\K[V_Q]$. From these data the
polynomials $v_{n-r+j}(Q,Z_{n-r+1})\ (2\le j\le r)$ can be easily
obtained by multiplication by $(\partial q/\partial
Z_{n-r+1})(Q,Z_{n-r+1})$ and modular reduction.

The polynomial $w_{n-r+1}(Q,Z_{n-r+1})$ can be computed as the
remainder of the product $(\partial{q}/\partial Z_{n-r+1})
(Q,Z_{n-r+1})\cdot{v}_{n-s+1}(Q,\!Z_{n-r+1})$ modulo
${q}(Q,\!Z_{n-r+1})$. Then, taking into account that the
identities $Y_{n-r+j}=\widehat{w}_{n-r+j}(1,Y_{n-r+1})$ and
$Y_{n-r+1}=v_{n-r+1}(Z_{n-r+1})$ hold in $\K[V_Q]$ for $2\le j\le
r$, we conclude that the polynomial $w_{n-r+j}(Q,Z_{n-r+1})$
equals the remainder of
$\widehat{w}_{n-r+j}\big(1,v_{n-r+1}(Z_{n-r+1})\big)$ modulo
$q(Q,Z_{n-r+1})$ for $2\le j\le r$. Therefore, the polynomials
$w_{n-r+j}(Q,Z_{n-r+1})$ $(2\le j\le r)$ can be computed with
$O(\delta_r\mathcal{U}(\delta_r))$ arithmetic operations in $\K$,
using at most $O(\delta_r^2)$ arithmetic registers.

Putting together the complexity and probability of success of each
step of the procedure above finishes the proof of the proposition.
\end{proof}
%
%
\section{The computation of a rational point of $V$}
\label{section:comput_rat_point}
In this section we exhibit a probabilistic algorithm which
computes a rational point of the variety $V:=V_r$. For this
purpose, let $\K$ be the finite field extension of $\fq$
introduced in Section \ref{section:computation_geo_sol} and assume
that we are given $\cfq$--independent linear forms $Z_1\klk
Z_{n-r+1},Y_{n-r+2}\klk Y_n\in\cfq[X]$, with $Z_1\klk
Z_{n-r+1}\in\fq[X]$ and $Y_{n-r+2}\klk Y_n\in\K[X]$, and a point
$Q:=(Q_1\klk Q_{n-r})\in\fq^{n-r}$, such that the linear projection
mapping $\pi:V\to\A^{n-r}$ determined by $Z_1\klk Z_{n-r+1}$
is a finite morphism and $Q$ is a lifting point of $\pi$. Furthermore,
assume that we are given
polynomials $q(Q,Z_{n-r+1})\in\fq[Z_{n-r+1}]$,
$v_{n-r+j}(Q,Z_{n-r+1})\in\K[Z_{n-r+1}]$ ($2\le j\le r$) which form a geometric
solution of the lifting fiber $V_Q$, as provided by Proposition
\ref{prop:fq_geo_sol}.

Let $\omega:=(\omega_1\klk \omega_{n-r})$ be an arbitrary point of
$\A^{n-r}$, let $L_\omega\subset\A^{n}$ be the
$(r+1)$--dimensional affine linear variety parametrized by
$Z_j=\omega_j T+Q_j$ $(1\le j\le n-r)$, $Z_{n-r+1}=Z_{n-r+1}$ and
$Y_{n-r+j}=Y_{n-r+j}$ $(2\le j\le r)$, and let
$\mathcal{C}_\omega:= V \cap L_\omega$. We may consider
$\mathcal{C}_\omega$ as the affine subvariety of $\A^{r+1}$
defined by the set of common zeros of the polynomials $F_j(\omega
T+Q,Z_{n-r+1}, Y_{n-r+2}\klk Y_n)$ $(1\le j\le r)$. With this
interpretation, let $\pi_\omega:V\to\A^1$ be the projection
mapping induced by $T$. We have the following result:
\begin{lemma}\label{lemma:prep_rat_point}
$\mathcal{C}_\omega$ is an equidimensional variety of $\A^{r+1}$
of dimension 1 and degree $\delta_r$, $\pi_\omega$ is a finite
morphism and 0 is an unramified value of $\pi_\omega$.
\end{lemma}
\begin{proof}
Lemma \ref{lemma:morph_bir} shows that the projection mapping
$\widetilde{\pi}:V\to \A^{n-r+1}$ defined by $Z_1\klk Z_{n-r+1}$
induces an isomorphism between a Zariski--dense open subset
$V\setminus \widetilde{V}$ of $V$ and a Zariski--dense open subset
$W\setminus \widetilde{W}$ of the hypersurface $W$ of $\A^{n-r+1}$
defined by $q(Z_1\klk Z_{n-r+1})$. Let
$\widetilde{L}_\omega\subset\A^{n-r+1}$ be the affine linear
variety parametrized by $Z_j=\omega_j T+Q_j$ $(1\le j\le n-r)$ and
$Z_{n-r+1}=Z_{n-r+1}$. Then we have that
$\widetilde{\pi}|_{L_\omega}:V\cap L_\omega\to
\widetilde{L}_\omega$ maps $\mathcal{C}_\omega=V\cap L_\omega$
into the hypersurface $W\cap \widetilde{L}_\omega$ of
$\widetilde{L}_\omega$ defined by the polynomial
$h(T,Z_{n-r+1}):=q(\omega T+Q,Z_{n-r+1})$. This shows that
$\widetilde{\pi}(\mathcal{C}_\omega)$ is an equidimensional
variety of dimension $1$. Furthermore,
$\widetilde{\pi}|_{L_\omega}$ maps $\widetilde{V}$ into
$\widetilde{W}\cap \widetilde{L}_\omega=V(\partial h/\partial
Z_{n-r+1})$. Taking into account that $h(0,Z_{n-r+1})=
q(Q,Z_{n-r+1})$ is squarefree, we easily conclude that
$W\cap\widetilde{W}\cap \widetilde{L}_\omega$ is a
zero--dimensional variety of $\widetilde{L}_\omega$. This shows
that $(W\setminus\widetilde{W})\cap \widetilde{L}_\omega$ is a
Zariski--dense open subset of $W\cap \widetilde{L}_\omega$ and
hence equidimensional of dimension 1. Therefore the Zariski--dense
open subset $(V\setminus\widetilde{V})\cap
{L}_\omega=\widetilde{\pi}^{-1}\big((W\setminus\widetilde{W})\cap
\widetilde{L}_\omega\big)$ of $\mathcal{C}_\omega$ is
equidimensional of dimension 1, which implies that
$\mathcal{C}_\omega$ is equidimensional of dimension 1.

The fact that the injective mapping $\cfq[Z_1\klk
Z_{n-r}]\hookrightarrow \cfq[V]$ induces an integral ring
extension implies that $\cfq[T]
\hookrightarrow\cfq[\mathcal{C}_\omega]$ is an injective mapping
which induces an integral extension ring, showing thus that
$\pi_\omega$ is a finite morphism. From the B\'ezout inequality
(\ref{equation:Bezout}), we see that $\deg \mathcal{C}_\omega\le
\delta_r$ holds. On the other hand, since $\pi_\omega^{-1}(0)=V_Q$
holds, we have $\delta_r=\deg V_Q\le\deg \mathcal{C}_\omega$. We
conclude that $\deg \mathcal{C}_\omega= \delta_r$ holds and 0 is
an unramified value of $\pi_\omega$.
\end{proof}

Let $\omega\in\fq^{n-r}$. Arguing as in the proof of Lemma
\ref{lemma:sucesionregular}, we easily conclude that $F_1(\omega
T+Q,Z_{n-r+1},Y_{n-r+2}\klk Y_n)\klk F_r(\omega
T+Q,Z_{n-r+1},Y_{n-r+2}\klk Y_n)$ form a regular sequence of
$\fq[T,Z_{n-r+1},Y_{n-r+2}\klk Y_n]$ and generate a radical ideal
of $\fq[T,Z_{n-r+1},Y_{n-r+2}\klk Y_n]$. Therefore, applying the
algorithm underlying Proposition \ref{prop:lifting_to_curve}, we
obtain elements $q(\omega T+Q,Z_{n-r+1})\in\fq[T,Z_{n-r+1}]$,
$v_{n-r+j}(\omega T+Q,Z_{n-r+1})\in\K[T,Z_{n-r+1}]$ $(2\le j\le
r)$, which define a geometric solution of the curve
$\mathcal{C}_\omega$.

Our intention is to find a rational point of the curve
$\mathcal{C}_\omega$. For this purpose, we are going to find a
rational point of the plane curve $W_\omega$ defined by the
polynomial $h:=q(\omega T+Q,Z_{n-r+1})$, which does not belong to
the plane curve $\widetilde{W}_\omega$ defined by the polynomial
$\partial h/\partial Z_{n-r+1}$. Let $\widetilde{\pi}_\omega:
\mathcal{C}_\omega\to \A^2$ be the mapping defined by
$T,Z_{n-r+1}$. From Lemma \ref{lemma:morph_bir} we deduce that
$\widetilde{\pi}_\omega$ induces a birational mapping
$\widetilde{\pi}_\omega: \mathcal{C}_\omega\to W_\omega$, whose
inverse is an $\fq$--definable rational mapping defined on
$W_\omega\setminus\widetilde{W}_\omega$, which can be easily
expressed in terms of the polynomials $v_{n-r+j}(\omega
T+Q,Z_{n-r+1})$ $(2\le j\le r)$. Therefore, using this inverse we
shall be able to obtain a rational point of our input variety $V$.
Unfortunately, the existence of a rational point of the plane
curve $W_\omega$ cannot be asserted if $W_\omega$ does not have at
least one absolutely irreducible component defined over $\fq$.

In order to assure that this condition holds, let
$C\in\cfq[\Omega_1 \klk \Omega_{n-r}]$ be the (nonzero) polynomial
of the statement of Theorem \ref{theorem:bertini}. Recall that $C$
has degree bounded by $2\delta_r^4$. Then Theorem
\ref{theorem:bertini} shows that, for any $\omega\in\fq^{n-r}$ with
$C(\omega)\not=0$, the curve $W_\omega$ is absolutely irreducible.

Assume as in Section \ref{subsec:fq_def_geo_sol} that
$q>8n^2d\delta_r^4$ holds. Theorem \ref{th:Zippel_Schwartz}
shows that a random choice of $\omega$ in $\fq^{n-r}$ satisfies
the condition $C(\omega)\not=0$ with probability at least
$1-1/72$. Then the (deterministic) algorithm underlying
Proposition \ref{prop:lifting_to_curve} yields a geometric
solution of the curve $\mathcal{C}_\omega$. We summarize the above
considerations in the following result:
\begin{proposition}\label{prop:prep_rat_point} Let $q>8n^2d
\delta_r^4$. There exists a probabilistic Turing machine $M$ which
has as input
\begin{itemize}
\item a \slp\ using space $\mathcal{S}$ and time $\mathcal{T}$
which represents the polynomials $F_1\klk F_r$,
\item the dense representation of elements of $\K[Z_{n-r+1}]$
which form a geometric solution of the lifting fiber $V_Q$,
as provided by Proposition \ref{prop:fq_geo_sol},
\end{itemize}
and outputs the dense representation of elements $q(\omega
T+Q,Z_{n-r+1})\in \fq[T,\!Z_{n-r+1}]$, $v_{n-r+j}(\omega
T+Q,Z_{n-r+1})\in\K[T,Z_{n-r+1}]$ $(2\le j\le r)$, which form a
geometric solution of the absolutely irreducible $\fq$--curve
$\mathcal{C}_\omega$. The Turing machine $M$ runs in space
$O\big((\mathcal{S}+n)\delta^2\log(q\delta)\big)$ and time
$O\big((n\mathcal{T}+n^5)\mathcal{U}(\delta)^2\mathcal{U}(\log
(q\delta))\big)$ and outputs the right result with probability at
least $1-1/72$.
\end{proposition}
%
%
\subsection{Computing a rational point of a plane curve}
\label{subsec:comput_point_curve}
In this subsection we exhibit a probabilistic algorithm which
computes a rational point of the curve $\mathcal{C}_\omega\subset
V $ previously defined.

Let $h:=q(\omega T+Q,Z_{n-r+1})$. Recall that $h$ is an absolutely
irreducible polynomial of $\fq[T,Z_{n-r+1}]$ of degree
$\delta_r>0$. Let as in the previous section $W_\omega,
\widetilde{W}_\omega\subset\A^2$ denote the plane curves defined
by $h$ and $\partial h/\partial Z_{n-r+1}$ respectively. As
remarked in the previous section, our aim is to compute a point in
the set $(W_\omega\setminus\widetilde{W}_\omega)\cap\fq^2$, from
which we shall immediately obtain a rational of point $V$.

The classical Weil's estimate on the number of rational points of
an absolutely irreducible projective plane curve \cite{Weil48}
implies that the set of rational points of $W_\omega$ satisfies
the estimate (see e.g. \cite{Schmidt74}):
$$|\#(W_\omega\cap\fq^2)-q|\le (\delta_r-1)(\delta_r-2)q^{1/2}+
\delta_r+1\leq \delta_r^2 q^{1/2}.$$
In particular, we deduce the lower bound $\#(W_\omega\cap\fq^2)\ge
q-\delta_r^2 q^{1/2}$.

On the other hand, by the absolute irreducibity of $h$ we conclude
that $h$ has no nontrivial common factor with $\partial h/\partial
Z_{n-r+1}$, which implies that $W_\omega\cap\widetilde{W}_\omega$
is a zero--dimensional variety. By the B\'ezout inequality we have
$\deg(W_\omega\cap\widetilde{W}_\omega)\le \delta_r(\delta_r-1)$,
which implies $\#(W_\omega\cap\widetilde{W}_\omega\cap\fq^2)\le
\delta_r(\delta_r-1)$. Combining this upper bound with our
previous lower bound, we conclude that the following estimate
holds:
\begin{equation}\label{eq:lower_bound_rat_points}
\#\big((W_\omega\setminus\widetilde{W}_\omega)\cap\fq^2\big)\ge
q-q^{1/2}\delta_r^2-\delta_r^2.\end{equation}
Assume that $q>8n^2d\delta_r^4$ holds. Then it is easy to see that
the right--hand side of (\ref{eq:lower_bound_rat_points}) is a
strictly positive real number, which implies that there exists at
least one rational point of
$W_\omega\setminus\widetilde{W}_\omega$.

Our purpose is to find a value $a\in\fq$ for which there exists a
rational point $(W_\omega\setminus \widetilde{W}_\omega)\cap\fq^2$
of the form $(a,z_{n-r+1})$. In order to find such value $a$, we
observe that for any $a\in\fq$ there exists at most $\delta_r$
points $(t,z_{n-r+1})\in W_\omega\setminus\widetilde{W}_\omega$
with $t=a$. Combining this observation with
(\ref{eq:lower_bound_rat_points}), we conclude that the following
estimate holds:
$$ \#\big\{a\in\fq:(W_\omega\setminus\widetilde{W}_\omega)
\cap\{T=a\}\cap \fq^2\not=\emptyset\big\}\geq
\frac{q-q^{1/2}\delta_r^2-\delta_r^2}{\delta_r}. $$
From this we immediately deduce the following lower bound on the
probability of finding at random a value $a$ for which there
exists a rational point with $t=a$:
\begin{equation}\label{eq:prob_estimate}
Prob\Big(a\in\fq:(W_\omega\setminus\widetilde{W}_\omega)
\cap\{T=a\}\cap \fq^2\not=\emptyset\Big)\geq\frac
{q-q^{1/2}\delta_r^2-\delta_r^2}{q\delta_r}.
\end{equation}
Let $q>8n^2d\delta_r^4$. Then the probability estimate
(\ref{eq:prob_estimate}) implies that, after at most $\delta_r$
random choices, we shall find a value $a\in \fq$ for which there
exists a rational point of $W_\omega\setminus
\widetilde{W}_\omega$ of the form $(a,z_{n-r+1})$ with probability
at least $1-2q^{-1/2}\delta_r^2\ge 1-1/6$. Having this value
$a\in\fq$, applying e.g. \cite[Corollary 14.16]{GaGe99} we see
that the computation of the value $z_{n-r+1}\in\fq$ can be reduced
to gcd computations and factorization in $\fq[Z_{n-r+1}]$. Our
next result describes the algorithm we have just outlined.
\begin{proposition}\label{prop:comput_rat_point} Let
$q>8n^2d\delta_r^4$. Then there exists a probabilistic Turing
machine $M$ which has as input a geometric solution of the plane
curve $\mathcal{C}_\omega$, as provided by Proposition
\ref{prop:prep_rat_point}, and outputs a rational point of
$\mathcal{C}_\omega$. The Turing machine $M$ runs in space
$O(\delta_r\log q\log(q\delta))$ and time
$O\big(n\delta_r\mathcal{U}(\delta_r)\log q\,\mathcal{U}(\log
(q\delta))\big)$, and outputs the right results with probability at
least $1-25/144$.
\end{proposition}
\begin{proof}
For $a\in \fq$, let
$h_a^*:=\gcd\big(h(a,Z_{n-r+1}),Z_{n-r+1}^q-Z_{n-r+1}\big)
\in\fq[Z_{n-r+1}]$. From \cite[Corollary 11.16]{GaGe99} we have
that the computation of $h_a^*$ can be performed with
$O\big(\mathcal{U}(\delta_r)\log q\big)$ operations in $\fq$,
storing $O(\delta_r\log q)$ elements of $\fq$. Furthermore,
deciding whether $h(a,Z_{n-r+1})$ is a squarefree polynomial
requires $O\big(\mathcal{U}(\delta_r)\big)$ operations in $\fq$,
storing $O(\delta_r)$ elements of $\fq$. From the
probability estimate (\ref{eq:prob_estimate}) we see that, after at
most $\delta_r$ random choices, with probability at least $1-1/6$
we shall find a value $a\in\fq$ such that $h(a,Z_{n-r+1})$ is
squarefree and $h_a^*$ is a nonconstant polynomial of
$\fq[Z_{n-r+1}]$. Therefore, computing such $a\in\fq$ and the
polynomial $h_a^*$ requires at most $O\big(\delta_r\mathcal{U}
(\delta_r)\log q\big)$ operations in $\fq$, storing
$O(\delta_r\log q)$ elements of $\fq$.

Observe that $h_a^*$ factors into linear factors in
$\fq[Z_{n-r+1}]$. Therefore, applying \cite[Theorem 14.9]{GaGe99}
we see that the factorization of $h_a^*$ in $\fq[Z_{n-r+1}]$
requires $O(\mathcal{U}(\delta_r)\log q)$ operations in $\fq$,
storing at most $O(\delta_r\log q)$, and outputs the right result
with probability at most $1-1/144$. Any root $b\in\fq$ of $h_a^*$
yields a rational point $(a,b)\in\fq^2$ of $W_\omega\setminus
\widetilde{W}_\omega$.

Specializing the parametrizations of $Y_{n-r+j}$ $(2\le j\le r)$
by the zeros of $q(\omega T+Q,Z_{n-r+1})$ into the values $T=a$ and
$Z_{n-r+1}=b$,
we obtain a rational point of $\mathcal{C}_\omega$ (observe that
our choice of $a$ assures that such specializations are
well--defined). This completes the proof of the proposition.
\end{proof}
Now we can describe the whole algorithm computing a rational point
of the input variety $V:=V_r$. First, we execute the algorithm
underlying Theorem \ref{th:geo_sol} in order to obtain a geometric
solution of the lifting fiber $\lfiberr$. Then we obtain a
geometric solution of the lifting fiber $V_Q$, and the absolutely
irreducible $\fq$--curve $\mathcal{C}_\omega$, applying the
algorithms underlying Propositions \ref{prop:fq_geo_sol} and
\ref{prop:prep_rat_point}. Finally, the algorithm of Proposition
\ref{prop:comput_rat_point} outputs a rational point of
$\mathcal{C}_\omega\subset V$. We summarize the result obtained in
the following corollary:
\begin{corollary}\label{coro:final} Let $q>8n^2d\delta_r^4$. Then
there exists a probabilistic Turing machine $M$, which takes as
input a \slp\ using space $\mathcal{S}$ and time $\mathcal{T}$
which represents the input polynomials $\fo{r}$, and outputs the
coordinates of a rational point of the variety $V:=V_r$. The
Turing machine $M$ runs in space
$O\big((\mathcal{S}+n+d)\delta\log q(\delta+\log (q\delta))\big)$
and time $O\big((n\mathcal{T}+n^5)\mathcal{U}(\delta)
\mathcal{U}(d\delta)\log q\,\mathcal{U}(\log (q\delta))\big)$, and
outputs the right result with probability at least $2/3>1/2$.
\end{corollary}
Let us remark that our algorithm can be easily extended to the
case of an equidimensional $\fq$--variety $V$ (given by a reduced
regular sequence), which has an absolutely irreducible component
defined over $\fq$. Indeed, the algorithm of Theorem
\ref{th:geo_sol} may be applied in this case, because it only
requires the variety $V$ to be equidimensional and given by a
reduced regular sequence. With a similar argument as in Theorem
\ref{theorem:bertini} and Proposition \ref{prop:prep_rat_point}, we
obtain a geometric solution of an $\fq$--curve $\mathcal{C}$,
contained in $V$, with at least one absolutely irreducible
component defined over $\fq$. Then, using fast algorithms for
bivariate factorization and absolute irreducibility testing (see
e.g. \cite{Kaltofen95a}), we compute such absolutely irreducible component,
to which we apply the algorithm underlying Proposition
\ref{prop:comput_rat_point}. Under the assumption that
$q>8n^2d\delta_r^4$ holds, the asymptotic complexity and
probability estimates of our algorithm in this case are the same
as in Corollary \ref{coro:final}.

\end{document}